\documentclass{amsart}

\usepackage{amssymb}
\usepackage{color}

 \newtheorem{theorem}{Theorem}[section]

 \newtheorem{proposition}[theorem]{Proposition}
 \newtheorem{lemma}[theorem]{Lemma}
 \newtheorem{corollary}[theorem]{Corollary}
 \newtheorem{remark}[theorem]{Remark}
 \newtheorem{definition}[theorem]{Definition}

 \newtheorem{problem}[theorem]{Problem}
 
 \newtheorem{condition}[theorem]{Condition}
 \newtheorem{assumption}[theorem]{Assumption}
 \numberwithin{equation}{section}
 \def\subrel#1#2{\mathrel{\mathop{#2}\limits_{#1}}}
 \def\qqq{\qed\vspace{2mm}}

\def\set{{\rm Set}}
\def\sset{{\rm sSet}}
\def\gdelta{G^{\delta}}
\def\Map{{\rm Map}}

\def\map{{\rm map}}

\def\sp{\Sigma{\rm Sp}}
\def\spk{\Sigma{\rm Sp}_{k}}
\def\isp{{\rm Sp}}
\def\ispk{{\rm Sp}_k}
\def\mod{{\rm Mod}}
\def\comod{{\rm Comod}}

\def\Ho{{\rm Ho}}

\begin{document}

\title
[
Discrete $G$-spectra and embeddings of module spectra]
{
Discrete $G$-Spectra
and embeddings of module spectra}
\author{Takeshi Torii}
\address{Department of Mathematics, 
Okayama University,
Okayama 700--8530, Japan}
\email{torii@math.okayama-u.ac.jp}

\thanks{The author would like to thank the anonymous referee for
many helpful comments and suggestions.
He also would like to thank Daniel Davis
for pointing out the mistake in Lemma~\ref{lemma:U-preserves-limit}
in the previous version of the paper.
This work was partially supported by 
JSPS KAKENHI Grant Number 25400092.}




\date{September 3, 2016\ ({\tt version~2.0})}

\begin{abstract}
In this paper we study the category
of discrete $G$-spectra for a profinite group $G$.
We consider an embedding of module objects 
in spectra into a category of 
module objects in discrete $G$-spectra,
and study the relationship between the embedding 
and the homotopy fixed points functor.
We also consider an embedding of module objects
in terms of quasi-categories, and 
show that the two formulations of embeddings
are equivalent in some circumstances.
\end{abstract}

\maketitle


\section{Introduction}

By the works of Morava~\cite{Morava},
Miller-Ravenel-Wilson~\cite{MRW},
Devinatz-Hopkins-Smith~\cite{DHS},
Hopkins-Smith~\cite{Hopkins-Smith},
Hovey-Strickland~\cite{Hovey-Strickland},
and many others,
the stable homotopy category 
is intimately related to the theory
of formal groups,
and the stable homotopy category localized at a prime $p$
has a filtration of full subcategories
corresponding to the height of formal groups.
The $n$th full subcategory for non-negative integer $n$
is realized as the $E(n)$-local category,
that is,
the Bousfield localization
with respect to the $n$th Johnson-Wilson theory $E(n)$ at $p$.
The $n$th subquotient of the filtration is equivalent
to the $K(n)$-local category, that is,
the Bousfield localization of the stable homotopy category
with respect to the $n$th Morava $K$-theory $K(n)$ at $p$,
and we can consider
that the fundamental building blocks of the stable homotopy category
are the $K(n)$-local categories
for various $n$ and $p$.   
Therefore,
studying the $K(n)$-local category has central
importance in stable homotopy theory.
A basic tool to study the $K(n)$-local category
is the $K(n)$-local $E_n$-based Adams spectral sequence,
where $E_n$ is the $n$th Morava $E$-theory at $p$.
The $E_2$-page of this spectral sequence 
for a $K(n)$-local $E(n-1)$-acyclic spectrum $X$ is 
described as the continuous cohomology of 
the extended Morava stabilizer group $\mathbb{G}_n$
with coefficients in the discrete twisted 
$E_{n*}$-$\mathbb{G}_n$-module $(E_n)_*(X)$.
This suggests that the derived category
of discrete (or continuous) twisted $E_{n*}$-$\mathbb{G}_n$-modules
may be intimately related to the $K(n)$-local category.   

An algebraic model of the $E(n)$-local category
was constructed by Franke~\cite{Franke} for 
sufficiently large primes $p$
compared to $n$,
after the pioneering work by Bousfield~\cite{Bousfield1}.
Although Franke's theory gives an equivalence of categories
between the derived category of periodic chain complexes
in $E(n)_*(E(n))$-comodules and the $E(n)$-local category,
it does not give an equivalence of model categories.
Therefore,
it does not give a model of the homotopy theory of $E(n)$-local spectra,
and we would like to have a model of the homotopy theory
of $E(n)$-local spectra, and of $K(n)$-local spectra.
In this paper we propose a model of the homotopy theory of $K(n)$-local
spectra in the model category of $K(n)$-local
$F_n$-modules in the discrete symmetric $\mathbb{G}_n$-spectra,
where $F_n$ is a discrete model of $E_n$ constructed
by Davis in \cite{Davis}
and upgraded to a commutative monoid 
object in the category of discrete symmetric $\mathbb{G}_n$-spectra
by Behrens-Davis in \cite{Behrens-Davis}.
For this purpose,
we shall formulate embeddings of module spectra
in general setting.

Let $G$ be a profinite group.
The theory of $G$-Galois extensions 
of structured ring spectra was introduced by
Rognes in \cite{Rognes},
and he gave an interpretation of the results
in \cite{Devinatz-Hopkins} in terms 
of $G$-Galois extensions.
The category of discrete $G$-spectra 
was introduced by Davis in \cite{Davis},
and $G$-Galois extensions were also 
studied by Behrens-Davis in \cite{Behrens-Davis}
in terms of discrete $G$-spectra.
The algebraic Galois theory is related 
to the descent theory,
which has been important 
in algebraic geometry, number theory, category theory
and homotopy theory.
The framework of homotopical descent theory was developed by
Hess~\cite{Hess} and Lurie~\cite{Lurie2}.
The formulation of embeddings we discuss in this paper 
is some kind of homotopical descent
related to $G$-Galois extensions of spectra.

Let $\sp(G)$ be the model category of discrete symmetric $G$-spectra
and let $\sp(G)_k$ be its left Bousfield localization 
with respect to a spectrum $k$. 
Suppose we have a map $A\to B$
of monoids in discrete symmetric $G$-spectra,
where $G$ acts on $A$ trivially.
There is a functor 
\[ {\rm Ex}: \mod_A(\spk)\to\mod_B(\sp(G)_k) \] 
from the category $\mod_A(\spk)$
of $A$-modules in $\spk$
to the category $\mod_B(\sp(G)_k)$
of $B$-modules in $\sp(G)_k$
by the extension of scalars.
This functor has a right adjoint, 
which we can regard as a fixed points functor $(-)^{G}$.
Furthermore,
the pair of functors is a $\sp$-Quillen adjunction.
We denote by
$\mathbb{L}{\rm Ex}$
the total left derived functor of ${\rm Ex}$ and
by $(-)^{hG}$ the total right derived functor
of the fixed points functor $(-)^{G}$.
The functor $(-)^{hG}$ can be regarded 
as the homotopy fixed points functor.
Let $\mathbf{T}$ be the full subcategory
of the homotopy category
$\Ho(\mod_A(\spk))$ consisting of $X$ such that
the unit map $X\to (\mathbb{L}{\rm Ex}(X))^{hG}$
is an equivalence.
We show that the restriction 
of $\mathbb{L}{\rm Ex}$ to $\mathbf{T}$
is fully faithful as an $\Ho(\sp)$-enriched functor 
(Proposition~\ref{prop:fundamental-proposition-module-embedding}).
Using 
the $K(n)$-local $\mathbb{G}_n$-Galois extension 
$L_{K(n)}S\to F_n$ formulated in \cite[\S8]{Behrens-Davis}
and the result in \cite{Davis-Torii},
we obtain the following theorem.

\begin{theorem}
[{Theorem~\ref{theorem:K(n)-local-embedding}}]
The total left derived functor 
\[ \mathbb{L}{\rm Ex}:\Ho(\sp_{K(n)})\longrightarrow
                \Ho(\mod_{F_n}(\sp(\mathbb{G}_n)_{K(n)})) \]
is fully faithful as an $\Ho(\sp)$-enriched functor.  
\end{theorem}

Next we consider embeddings of modules in quasi-categories.
Model categories are models of homotopy theories
and contain rich homotopy theoretic structures.
But they are sometimes rigid and hard to work with
since they contain auxiliary structures
such as cofibrations and fibrations.
Simplicial categories and topological categories are
also models of homotopy theories.
But they are also sometimes rigid and hard to work with
since they have a strict associative composition law
of mapping spaces.
Quasi-categories are yet other models of homotopy theories.
Quasi-categories were introduced 
by Boardman-Vogt~\cite{Boardman-Vogt}
and developed by Joyal~\cite{Joyal,Joyal2} and 
Lurie~\cite{Lurie,Lurie2}.
The definition of quasi-categories is a little strange
at first glance but
actually quasi-categories are flexible, easy to work with,
and well-developed.  
Therefore, to formulate embeddings of modules 
in terms of quasi-categories is important and will be useful for 
later applications. 

Let $\ispk$ be the underlying quasi-category of 
the simplicial model category $\spk$.
Suppose we have a map $A\to E$ 
of algebra objects in $\ispk$.
We have a functor from the quasi-category
$\mod_A(\ispk)$ of $A$-modules 
to the quasi-category $\mod_E(\ispk)$ of $E$-modules
by the extension of scalars. 
This defines a comonad $\Theta$ on the $E$-modules, and 
we can consider the quasi-category of comodules 
$\comod_{(E,\Theta)}(\ispk)$ 
over the comonad $\Theta$.
The extension of scalars functor factors through
$\comod_{(E,\Theta)}(\ispk)$,
and we obtain a functor
\[ {\rm Coex}: \mod_A(\ispk)\to \comod_{(E,\Theta)}(\ispk). \]
This functor has a right adjoint $P$,
which is a homotopical analogue
of the functor taking primitive elements.
Let $\mathbf{T}$ be the full subcategory
of $\mod_A(\ispk)$ consisting of $X$
such that the unit map $X\to P{\rm Coex}(X)$ 
is an equivalence.
We show that the restriction 
of the functor ${\rm Coex}$ to $\mathbf{T}$
is fully faithful 
(Proposition~\ref{prop:quasi-category-module-embedding}).
This result is related to the 
(effective) homotopic descent considered by Hess
in \cite[Def.~5.1]{Hess},
but the author is not sure if 
this result can be formulated in terms of
simplicial categories.
There is a map of quasi-categories
from the simplicial nerve of the simplicial category
of coalgebras in ${\rm Mod}_E(\sp_k)$
associated to the adjunction
$\mod_A(\sp_k)\rightleftarrows {\rm Mod}_E(\spk)$ 
to the quasi-category of comodules 
${\rm Comod}_{(E,\Theta)}(\ispk)$,
but it seems that this map may not be an equivalence
in general
since a comodule in ${\rm Comod}_{(E,\Theta)}(\ispk)$
which satisfies the comodule structure
up to higher coherent homotopies may not
be rectified as a coalgebra in $\mod_A(\sp_k)$
which satisfies the comodule structure on the nose.

Finally, we compare the two formulations 
of embeddings in terms of model categories and
of quasi-categories.
For this purpose we describe the underlying
quasi-categories of modules and algebras
in discrete symmetric $G$-spectra for a
profinite group $G$ in terms of
the quasi-categories of modules and algebras
in non-equivariant symmetric spectra.

Let $A\to B$ be a map of monoids in
discrete symmetric $G$-spectra,
where $G$ acts on $A$ trivially.
Let $U: \sp(G)_k\to \sp_k$ be the forgetful functor,
which has the right adjoint $V={\rm Map}_c(G,-)$
(see \S\ref{subsection: discrete-G-spectra} 
for the functor ${\rm Map}_c(G,-)$).
We have a map
\[ \Psi_{UB}: B(UB,A,UB)\longrightarrow
         U{\rm Map}_c(G,UB),\]
where $B(UB,A,UB)$ is the two-sided bar construction
of the $A$-module $UB$ 
(see \S\ref{subsection:equivalence-two-formulations}
for the construction of the map $\Psi_{UB}$).
Let $\mod_B(\isp(G)_k)$ be the underlying
quasi-category of the simplicial model category 
$\mod_B(\sp(G)_k)$.
We show that $\mod_B(\isp(G)_k)$
can be written as a quasi-category of
comodules.

\begin{theorem}[Corollary~\ref{corollary:mod-G-equivalent-comodules}]
\label{theorem:main-theorem-2}
Let $G$ be a profinite group that has finite virtual 
cohomological dimension.
We assume that the localization functor
$L_k$ is given as a composite $L_ML_T$ of
two localization functors $L_M$ and $L_T$,
where $L_T$ is a smashing localization and 
$L_M$ is a localization with respect to a finite spectrum $M$. 
If $\Psi_{UB}$ is a $k$-local equivalence,
then there is an equivalence of quasi-categories
\[ \mod_B(\isp(G)_k)\simeq
   \comod_{(UB,\Theta)}(\ispk).\]
\end{theorem}

As a corollary,
we obtain that the functor ${\rm Ex}$
is equivalent to the functor ${\rm Coex}$
under the equivalence between
$\mod_B(\isp(G)_k)$ and
$\comod_{(UB,\Theta)}(\ispk)$
(Corollary~\ref{corollary:equivalence-formulations}),
where we regard ${\rm Ex}$ as a functor of
the underlying quasi-categories.
This shows that the two formulations 
of embeddings
are equivalent.
In particular,
we show that the two formulations
are equivalent 
if $A\to B$ is a $k$-local $G$-Galois extension 
(Theorem~\ref{thm:Galis-mod-G-equivalent-comodules})
under the assumptions of Theorem~\ref{theorem:main-theorem-2}.

The organization of this paper is as follows:
In \S\ref{section:model_structure_spG}
we discuss the model structure on the category
of discrete symmetric $G$-spectra.
We show that the category of discrete symmetric $G$-spectra
is a proper, combinatorial, symmetric monoidal 
$\sp$-model category satisfying the monoid axiom.
We also discuss the Bousfield localization 
with respect to a spectrum with trivial $G$-action.
In \S\ref{section:embedding_model_module}
we discuss embeddings of modules 
into the category of 
discrete symmetric $G$-spectra. 
We also discuss the relationship between
the embeddings and the homotopy fixed points functors. 
In \S\ref{section:embedding-quasi-categories}
we consider embeddings of quasi-categories of modules
in spectra. 
For an adjunction of quasi-categories,
we can consider the quasi-category of comodules
over the comonad associated to the adjunction.
We show that some full subcategory
can be embed into the quasi-category of comodules.
In \S\ref{section:quasi-category-spG}
we study the quasi-category of discrete $G$-spectra.
We show that the quasi-category
of discrete $G$-spectra can be described as
a quasi-category of comodules.
Finally we show that the two formulations
of embeddings of module categories 
are equivalent in some circumstances.
In \S\ref{section:Galois-descent}
we discuss embeddings associated
to profinite $G$-Galois extensions.
We show that the two formulations are
equivalent for
profinite $G$-Galois extensions
under some conditions.

\section{Notation}

For a model category $\mathbf{M}$,
we denote by $\Ho(\mathbf{M})$
the homotopy category of $\mathbf{M}$.
For objects $X,Y\in\mathbf{M}$,
we denote by $[X,Y]_{\mathbf{M}}$
the set of morphisms in $\Ho(\mathbf{M})$.
For a simplicial model category $\mathbf{N}$,
we denote by
${\rm Map}_{\mathbf{N}}(X,Y)$
the mapping space (simplicial set) for $X,Y\in\mathbf{N}$ 
(see, for example,
\cite{Goerss-Jardine,Hirschhorn,Hovey}
for these concepts).
We denote by $\mathbf{N}^{\circ}$
the full simplicial subcategory of $\mathbf{N}$
consisting of objects that
are both fibrant and cofibrant
as in \cite{Lurie}.
The underlying quasi-category of 
$\mathbf{N}$
is defined to be $N(\mathbf{N}^{\circ})$,
where $N(-)$ is the simplicial nerve functor
(see \cite[1.1.5]{Lurie}
for the simplicial nerve functor).

We denote by $\sp$ 
the category of symmetric spectra
constructed by Hovey-Shipley-Smith in
\cite{HSS}.
We give $\sp$ the stable model structure
(see \cite[\S3]{HSS} for the definition
of the stable model structure on $\sp$).
We denote by $\isp$ the underlying quasi-category of $\sp$.
For a spectrum $k$,
we denote by $\spk$ 
the left Bousfield localization of $\sp$
with respect to $k$
(see, for example, \cite[Ch.~3]{Hirschhorn}
for the definition of left Bousfield localization
of a model category), 
and by $\ispk$ the underlying quasi-category
of $\spk$.
We denote by $\mathcal{S}$
the quasi-category of spaces,
which is the underlying quasi-category
of the category of simplicial sets
with the Kan model structure
(see \cite[1.2.16]{Lurie}).
For a quasi-category $\mathcal{C}$,
we have a mapping space (simplicial set)
${\rm Map}_{\mathcal{C}}(X,Y)$ for $X,Y\in \mathcal{C}$,
which is well-defined up to weak homotopy equivalence
(see \cite[1.2.2]{Lurie}).

\section{Model structure on the 
category of discrete symmetric $G$-spectra}
\label{section:model_structure_spG}

Let $G$ be a profinite group.
In this section we discuss model structure
on the category of discrete symmetric $G$-spectra.
We also study the Bousfield localization with
respect to a spectrum with trivial $G$-action.

\subsection{Discrete symmetric $G$-spectra}
\label{subsection: discrete-G-spectra}

Let $G$ be a profinite group.
In this subsection
we recall the model structure
on the category of discrete symmetric $G$-spectra
and study its properties.
We show that the category of discrete symmetric $G$-spectra
is a proper combinatorial 
symmetric monoidal $\sp$-model category
satisfying the monoid axiom.
We also compare the model category of discrete symmetric $G$-spectra
with that of non-equivariant symmetric spectra.

First, we recall the definition of a discrete symmetric $G$-spectrum
(see \cite[\S2.3]{Behrens-Davis}). 
We denote by $\set(G)$ the category of discrete $G$-sets.
A simplicial discrete $G$-set is a simplicial object 
in $\set(G)$.
The model structure on the category of simplicial discrete
$G$-sets was studied in \cite{Goerss}.
We denote by $\sset(G)_*$ the category of pointed 
simplicial discrete $G$-sets.
Let $\sset(G)_*^{\Sigma}$ be
the category of symmetric sequences in $\sset(G)_*$.
We can give a closed symmetric monoidal
structure on $\sset(G)_*^{\Sigma}$.
Let $S$ be the symmetric sequence given by
$S=(S^0,S^1,S^2,\ldots)$,
where $S^n$ is the $n$-sphere with trivial $G$-action.
The symmetric sequence $S$ is a commutative monoid
object in $\sset(G)_*^{\Sigma}$.
A discrete symmetric $G$-spectrum
is a module object in $\sset(G)_*^{\Sigma}$
over the commutative monoid $S$.
A map of discrete symmetric $G$-spectra is a map
of module objects.
We denote by $\sp(G)$ the category of discrete symmetric
$G$-spectra.

The category $\sp(G)$ is a complete, cocomplete, closed 
symmetric monoidal category with
$S$ as the unit object.  
We denote the monoidal structure by
$X\wedge Y = X\otimes_S Y$.
We have an adjoint pair 
\[ {\rm triv}: \sp\rightleftarrows 
      \sp(G): (-)^G,\]
where $(-)^G$ is the $G$-fixed points functor and
the functor ${\rm triv}(-)$ associates to 
a symmetric spectrum $X$
the discrete symmetric $G$-spectrum $X$ 
with trivial $G$-action.
Notice that the functor ${\rm triv}$ 
is a strong symmetric monoidal functor
and the functor $(-)^G$ is a lax symmetric monoidal functor.

Now we recall the model structure on $\sp(G)$
defined in \cite[\S2.3]{Behrens-Davis}.
Let $U: \sp(G)\to\sp$ be the forgetful functor.
A map $f:X\to Y$ in $\sp(G)$
is said to be 
\begin{itemize}
\item
a cofibration if $U(f)$ is a cofibration of symmetric spectra,
\item
a weak equivalence if $U(f)$ is a stable equivalence
of symmetric spectra, and 
\item
a fibration if it has
the right lifting property with respect to
all maps which are both cofibrations and
weak equivalences.
\end{itemize}  
With these definitions,
$\sp(G)$ is a left proper cellular model category
by \cite[Thm.~2.3.2]{Behrens-Davis}.
Recall that a model category is called left proper if 
every pushout of a weak equivalence along a cofibration
is a weak equivalence
(see \cite[Ch.~13]{Hirschhorn}).
A cellular model category is a cofibrantly
generated model category 
for which there are a set $I$ of generating
cofibrations and a set $J$ of generating trivial 
cofibrations such that
\begin{itemize}
\item
both the domains and the codomains of the elements
of $I$ are compact,
\item
the domains of the elements of $J$ are small relative
to $I$, and
\item
the cofibrations are effective monomorphisms
\end{itemize}
(see \cite{Hirschhorn} for these concepts,
in particular, \cite[Ch.~12]{Hirschhorn} for cellular model categories).


We shall show that $\sp(G)$ is a proper combinatorial model category.
Recall that a model category is called 
right proper if
every pullback of a weak equivalence along a fibration
is a weak equivalence, and
proper if it is both left proper and right proper
(see \cite[Ch.~13]{Hirschhorn}).
A model category 
is said to be combinatorial
if it is cofibrantly generated as a model category 
and is locally presentable as a category
(see, for example, \cite[\S2]{Dugger1} or
\cite[A.2.6]{Lurie} for combinatorial model categories).
Since $\sp(G)$ is a cofibrantly generated model category, 
it suffices to show that $\sp(G)$ is locally presentable
in order to show that $\sp(G)$ is combinatorial.
Recall that a category is locally $\lambda$-presentable
for a regular cardinal $\lambda$
if it is cocomplete and has a set $C$
of $\lambda$-compact objects such that
every object is a $\lambda$-filtered colimit of 
objects in $C$.
A category is called locally presentable if
it is locally $\lambda$-presentable for
some regular cardinal $\lambda$
(see \cite[Ch.~1]{Adamek-Rosicky}
for locally presentable categories).

\begin{theorem}
\label{thm:spG-combinatorial}
The category $\sp(G)$ is locally presentable.
Hence $\sp(G)$ is a combinatorial model category.
\end{theorem}

\proof
The category ${\rm Set}(G)$ 
of discrete $G$-sets is locally $\aleph_0$-presentable,
where $\aleph_0$ is the first infinite cardinal.
A discrete $G$-set is $\aleph_0$-compact
if and only if the underlying set is finite.
The full subcategory of finite discrete $G$-sets
is essentially small.
We denote by $\mathcal{A}$ the opposite
category of a skeleton of the full subcategory
of finite discrete $G$-sets.
For a small category $\mathcal{C}$
and a category $\mathcal{D}$,
we denote by ${\rm Fun}(\mathcal{C},\mathcal{D})$
the functor category and by
${\rm Fun}^{\rm Lex}(\mathcal{C},\mathcal{D})$
the full subcategory of finite-limit preserving functors.
By \cite[Thm.~1.46 and its proof]{Adamek-Rosicky},
we see that the Yoneda map
${\rm Set}(G)\to {\rm Fun}(\mathcal{A},{\rm Set})$
given by $X\mapsto {\rm Hom}_{{\rm Set}(G)}(-,X)$
for $X\in {\rm Set}(G)$
induces an equivalence of categories
between ${\rm Set}(G)$ and 
${\rm Fun}^{\rm Lex}(\mathcal{A},{\rm Set})$.

For a (discrete) symmetric ($G$-)spectrum $Y$,
we denote by $Y_{k,l}$
the set of $l$-simplexes of the $k$th simplicial
set $Y_k$ of $Y$.
We can define a functor
$F: \Sigma{\rm Sp}(G)\to
   {\rm Fun}^{\rm Lex}(\mathcal{A},\Sigma{\rm Sp})$
by 
\[ F(X)(A)_{k,l}={\rm Hom}_{{\rm Set}(G)}(A,X_{k,l}) \]
for $A\in\mathcal{A}$ and $X\in \Sigma{\rm Sp}(G)$
with obvious structure maps.
We can verify that the functor $F$ is an equivalence of
categories by using the equivalence
${\rm Set}(G)\stackrel{\simeq}{\to} 
{\rm Fun}^{\rm Lex}(\mathcal{A},{\rm Set})$.

Recall that the category $\Sigma{\rm Sp}$
of symmetric spectra is locally presentable.
This follows, for example, from
\cite[1.2.10, 3.2.13, and 5.1.6]{HSS}.
See also \cite[p.474]{Shipley}.
The theorem follows from the fact that
for any small category $\mathcal{C}$
and any locally presentable category $\mathcal{D}$,
the category ${\rm Fun}^{\rm Lex}(\mathcal{C},\mathcal{D})$
is locally presentable
by \cite[1.53 and 1.50(1)]{Adamek-Rosicky}.
\qqq

\begin{proposition}
The model category $\sp(G)$ is proper.
\end{proposition}

\proof
It suffices to show that $\sp(G)$ is right proper.
Since $\sp$ is right proper
by \cite[Thm.~5.5.2]{HSS},
this follows from the fact that
the forgetful functor $U$ preserves
fiber products and detects weak equivalences.
\qqq

Next we consider the compatibility of the monoidal structure 
and the model structure on $\sp(G)$. 
By the definition of the model structure on $\sp(G)$
and the fact that the composition
$U\circ {\rm triv}$ is the identity functor,
the functor ${\rm triv}$ preserves cofibrations and weak equivalences,
and hence the adjoint pair of functors
\[ {\rm triv}: \sp\rightleftarrows \sp(G): (-)^G\]
is a Quillen adjunction.

We shall recall the definition of a symmetric monoidal 
Quillen adjunction (see \cite[\S4.2]{Hovey})
and show that the pair 
$({\rm triv},(-)^G)$
is a symmetric monoidal Quillen adjunction.
Let $\mathbf{M}$ and $\mathbf{N}$ be 
a symmetric monoidal model categories.
A Quillen adjunction 
$F:\mathbf{M}\rightleftarrows \mathbf{N}:G$
is said to be 
a symmetric monoidal Quillen adjunction
if the left Quillen functor $F$ is strong symmetric monoidal
and the map $F(q): F(Q\mathbb{I})\to F(\mathbb{I})$ 
is a weak equivalence,
where $\mathbb{I}$ is a unit object in $\mathbf{M}$
and $q: Q\mathbb{I}\to \mathbb{I}$ 
is a cofibrant replacement of $\mathbb{I}$.  
We say that the left adjoint of a symmetric monoidal
Quillen adjunction is a symmetric monoidal left Quillen functor.
Since the functor ${\rm triv}$ is strong symmetric monoidal
and the sphere spectrum is cofibrant in $\sp$,
we see that the pair $({\rm triv}, (-)^G)$
is a symmetric monoidal Quillen adjunction.

Let $\mathbf{C}$ be a symmetric monoidal model category.
We also recall the definitions of a symmetric monoidal 
$\mathbf{C}$-model category
and a symmetric monoidal $\mathbf{C}$-Quillen adjunction
(see \cite[\S4.2]{Hovey}).
A model category $\mathbf{M}$ 
is said to be  
a symmetric monoidal $\mathbf{C}$-model category 
if it is a symmetric monoidal model category
equipped with a strong symmetric monoidal left Quillen functor 
$i: \mathbf{C}\to \mathbf{M}$. 
Let $\mathbf{M}_1$ and $\mathbf{M}_2$ be
symmetric monoidal $\mathbf{C}$-model categories
equipped with symmetric monoidal left Quillen functors
$i_1:\mathbf{C}\to\mathbf{M}_1$ and 
$i_2:\mathbf{C}\to\mathbf{M}_2$, respectively.
A symmetric monoidal $\mathbf{C}$-Quillen
adjunction between $\mathbf{M}_1$ and $\mathbf{M}_2$
is a symmetric monoidal Quillen adjunction
$F:\mathbf{M}_1\rightleftarrows \mathbf{M}_2:G$
together with a symmetric monoidal 
natural isomorphism between $F\circ i_1$ and $i_2$.
We say that the left adjoint of a symmetric monoidal
$\mathbf{C}$-Quillen adjunction 
is a symmetric monoidal left $\mathbf{C}$-Quillen functor.

\begin{theorem}
\label{thm:weak-stable-symmetric-monoidal}
The category 
$\sp(G)$ is a symmetric monoidal 
$\sp$-model category.
The adjoint pair $({\rm triv},(-)^G)$
is a symmetric monoidal $\sp$-Quillen adjunction.
\end{theorem}

\proof
The model category $\sp(G)$
is a symmetric monoidal model category 
by \cite[Thm.~8.11]{Hovey2}.
The theorem follows from the fact that 
${\rm triv}: \sp\to \sp(G)$ is a symmetric monoidal 
left Quillen functor.
\qed


\bigskip

Now we shall verify that
the symmetric monoidal model category
$\sp(G)$ satisfies the monoid axiom.
We recall the monoid axiom on a monoidal model category
(see, for example, \cite[Def.~3.3]{Schwede-Shipley}).
Let $\mathbf{M}$ be a monoidal model category
with tensor product $\wedge$.
For a class $\mathcal{I}$ of maps in $\mathbf{M}$,
we denote by $\mathcal{I}\wedge \mathbf{M}$
the class of maps of the form
\[ A\wedge Z\longrightarrow B\wedge Z \]
for $A\to B$ a map in $\mathcal{I}$ and $Z$ an object of $\mathbf{M}$.
For a class $\mathcal{J}$ of maps in $\mathbf{M}$,
we denote by $\mathcal{J}\mbox{\rm -cof}_{\rm reg}$ the class
of maps obtained from the maps of $\mathcal{J}$ by cobase change
and transfinite composition
(see, for example, 
\cite[10.2]{Hirschhorn}
or \cite[2.1.1]{Hovey}
for the definition of a
transfinite composition).  
We say that a monoidal model category $\mathbf{M}$ satisfies
the monoid axiom if every map in
\[ (\{\mbox{\rm trivial cofibrations}\}\wedge \mathbf{M})
   \mbox{\rm -cof}_{\rm reg} \]
is a weak equivalence,
where $\{\mbox{\rm trivial cofibrations}\}$ is
the class of trivial cofibrations.

\begin{proposition}\label{prop:weak-stable-monoid-axiom}
The symmetric monoidal model category
$\sp(G)$ satisfies the monoid axiom. 
\end{proposition}

\proof
We have to show that every map in
\[ (\{\mbox{\rm trivial cofibrations}\}\wedge 
     \sp(G))\mbox{\rm -cof}_{\rm reg} \]
is a weak equivalence.
By the definition of the model structure on $\sp(G)$,
it suffices to show that the underlying map
is a stable equivalence of symmetric spectra.
Since the underlying map of a trivial cofibration 
is a trivial cofibration of symmetric spectra,
the proposition follows from the fact that
the category of symmetric spectra 
with stable model structure
satisfies the monoid axiom
by \cite[Thm.~5.4.1]{HSS}.
\qed

\bigskip

Now we shall compare the category of
discrete symmetric $G$-spectra
with that of non-equivariant symmetric spectra.
For this purpose,
we introduce
a (non-discrete) symmetric $G$-spectra.
We denote by $\gdelta$ the group $G$ with discrete topology.
Let $\sp(\gdelta)$ be the category of symmetric spectra
with (continuous) $G^{\delta}$-action.
The continuous homomorphism $\gdelta\to G$
induces a functor
$(-)^{\delta}: \sp(G)\to \sp(\gdelta)$. 
For $X\in\sp(\gdelta)$,
we denote by $dX$ the largest discrete $G$-subspectrum of $X$,
that is, 
\[ dX={\rm colim}_H X^H, \]
where $H$ ranges over all open subgroups $H$ of $G$.
We can regard $d$ as a functor 
$d: \sp(\gdelta)\to\sp(G)$.
Notice that we have an adjoint pair
\[ (-)^{\delta}: \sp(G)\rightleftarrows\sp(\gdelta):d.\]

Let $U:\sp(G)\to\sp$
be the forgetful functor.
The functor $U$ is strong symmetric monoidal 
and has a right adjoint 
\[ V: \sp\to\sp(G).\]
We shall explicitly describe the right adjoint $V$.
For a symmetric spectrum $Y$ (with $G$-action),
we denote by $Y_{k,l}$
the set of $l$-simplexes of the $k$th simplicial set $Y_k$ of $Y$.  
For a set $A$, we denote by
${\rm Map}(G,A)$ the $G$-set of all maps
from $G$ to $A$ with $G$-action
given by $(g\cdot \theta) (g')=
\theta(g'g)$ for $g,g'\in G$ and
$\theta\in {\rm Map}(G,A)$.
For a symmetric spectrum $X$,
we can define an object ${\rm Map}(G,X)$
of $\sp(G^{\delta})$ by
${\rm Map}(G,X)_{k,l}={\rm Map}(G,X_{k,l})$
with obvious structure maps.
We define the discrete symmetric $G$-spectrum 
${\rm Map}_c(G,X)$ by
\[ {\rm Map}_c(G,X)=d({\rm Map}(G,X)).\]
The right adjoint $V$ of the forgetful 
functor $U: \sp(G)\to \sp$ is given by 
\[ V(X)= \Map_c(G,X).\]

We can easily verify the following proposition.

\begin{proposition}
\label{prop:U-V-adjunction}
The adjoint pair of functors 
\[ U: \sp(G)\rightleftarrows\sp:V.\]
is a symmetric monoidal $\sp$-Quillen adjunction.
\end{proposition}

\proof
By the definition of the model structure 
on $\sp(G)$,
the forgetful functor $U$ preserves cofibrations and
weak equivalences, and hence
the adjoint pair $(U,V)$ of functors
is Quillen adjunction.
Since $U$ is a strong symmetric monoidal
and the sphere spectrum is cofibrant
in $\sp(G)$,
we see that the pair $(U,V)$ is a symmetric
monoidal Quillen adjunction.
Since the composition $U\circ {\rm triv}$
is the identity functor on $\sp$,
we see that the forgetful functor $U$ is 
a symmetric monoidal left $\sp$-Quillen functor.
This completes the proof.
\qqq


\subsection{Bousfield localization of $\sp(G)$}

Let $k$ be a symmetric spectrum.
In this subsection we study the left Bousfield localization
$\sp(G)_k$ of $\sp(G)$ with respect to $k$.
We show that $\sp(G)_k$
is a left proper combinatorial symmetric
monoidal $\sp$-model category satisfying the monoid axiom.
We also compare the model category $\sp(G)_k$
with the left Bousfield localization $\spk$
of non-equivariant symmetric spectra with respect to $k$.

First, we recall that 
there exists the left Bousfield localization on $\sp(G)$ 
with respect to $k$.
This follows, for example, from the fact that 
$\sp(G)$ is a left proper cellular model category.
We say that a morphism $f$ in $\sp(G)$
is a $k$-local equivalence if 
$U(f)$ is a $k$-local equivalence in $\sp$.
Let $W_k$ be the class of $k$-local equivalences in $\sp(G)$.
As in \cite[p.5015]{Behrens-Davis},
there exists the left Bousfield localization
$\sp(G)_k$ with respect to $W_k$,
and $\sp(G)_k$ is left proper and cellular.

\begin{proposition}
\label{prop:spGk-combinatorial}
The left Bousfield localization $\sp(G)_k$
is a left proper combinatorial simplicial model category.
\end{proposition}

\proof
By Theorems~\ref{thm:spG-combinatorial}
and \ref{thm:weak-stable-symmetric-monoidal},
and \cite[Thm.~2.3.2]{Behrens-Davis},
we see that $\sp(G)$ is a left proper 
combinatorial simplicial model category.
As in \cite[p.5015]{Behrens-Davis},
the class of $k$-local equivalences are that
of $f$-local equivalences for some map $f$.
Hence the left Bousfield localization
$\sp(G)_k$ is a left proper combinatorial simplicial
model category
by \cite[Prop.~A.3.7.3]{Lurie}.
\qqq

In the following of this paper we assume that
$k$ is cofibrant for simplicity.

\begin{theorem}
\label{thm:spgk-symmetric-monoidal-model-category}
The model category $\sp(G)_k$
is a symmetric monoidal $\sp$-model category.
\end{theorem}

\proof
First, we show that $\sp(G)_k$ is 
a symmetric monoidal model category.
Let $A\to B$ be a cofibration and let $X\to Y$
be a trivial cofibration in $\sp(G)_k$.
Since $k$ is cofibrant,
we see that $X\wedge k \to Y\wedge k$
is a trivial cofibration in $\sp(G)$.
By Theorem~\ref{thm:weak-stable-symmetric-monoidal},
$\sp(G)$ is a monoidal model category,
and hence the map
\[ (A\wedge Y\wedge k)\subrel{(A\wedge X\wedge k)}{\coprod}
 (B\wedge X\wedge k)\to B\wedge Y\wedge k \]
is a trivial cofibration.
We have an isomorphism
\[ \left((A\wedge Y)\subrel{(A\wedge X)}{\coprod}(B\wedge
   X)\right)\wedge k \cong
   (A\wedge Y\wedge k)\subrel{(A\wedge X\wedge k)}{\coprod}  
   (B\wedge X\wedge k), \]
and hence the map 
\[ (A\wedge Y)\subrel{(A\wedge X)}{\coprod}(B\wedge X)
  \to B\wedge Y\]
is a trivial cofibration
in $\sp(G)_k$.
Therefore, we see that $\sp(G)_k$ is 
a symmetric monoidal model category.

Next, we show that $\sp(G)_k$ is a symmetric
monoidal $\sp$-model category.
Recall that the functor 
${\rm triv}: \sp\to \sp(G)$ is 
a symmetric monoidal left Quillen functor. 
Furthermore,
the identity functor ${\rm id}: \sp(G)\to\sp(G)_k$
is a symmetric monoidal left Quillen functor.
Hence the composition ${\rm id}\circ{\rm triv}: \sp\to\sp(G)_k$ is 
also a symmetric monoidal left Quillen functor.
This shows that
$\sp(G)_k$ is a symmetric monoidal $\sp$-model category.
\qed

\begin{proposition}
\label{prop:spgk-satisfies-monoid-axiom}
The symmetric monoidal model category $\sp(G)_k$
satisfies the monoid axiom.
\end{proposition}

\proof
Recall that $W_k$ is the class of $k$-local equivalences. 
We let $C$ be the class of cofibrations in $\sp(G)_k$.
Since the functor $(-)\wedge k$ 
preserves all colimits,
we see that
\[ \left(((C\cap W_k)\wedge \sp(G))\mbox{\rm -cof}_{\rm reg}\right)
    \wedge k\subset 
    ((C\cap W)\wedge \sp(G))\mbox{\rm -cof}_{\rm reg},
   \]
where $W$ is the class of weak equivalences in $\sp(G)$.
By Proposition~\ref{prop:weak-stable-monoid-axiom},
$\sp(G)$ satisfies the monoid axiom.
Hence 
$((C\cap W)\wedge \sp(G))\mbox{\rm -cof}_{\rm reg}\subset W$.
This shows that
\[ ((C\cap W_k)\wedge \sp(G))\mbox{\rm -cof}_{\rm reg}
   \subset W_k.\]
This completes the proof.
\qed

\bigskip

Now we compare the model category $\sp(G)_k$
with the left Bousfield localization $\spk$
of non-equivariant symmetric spectra with respect to $k$.

\begin{proposition}
\label{prop:(i,(-)^G-adjunction-localizaed-version)}
The adjoint pair of functors
\[ {\rm triv}: \spk\rightleftarrows \sp(G)_k: (-)^G\]
is a symmetric monoidal $\sp$-Quillen adjunction.
\end{proposition}

\proof
It suffices to show that ${\rm triv}$ preserves $k$-local equivalences.
Let $f$ be a $k$-local equivalence in $\sp$.
Since $U({\rm triv}(f))=f$,
we see ${\rm triv}(f)\in W_k$.
Hence ${\rm triv}$ preserves $k$-local equivalences.
\qed

\begin{proposition}
\label{proposition:k-local-U-V-adjunction}
The adjoint pair of functors 
\[ U: \sp(G)_k\rightleftarrows\spk:V.\]
is a symmetric monoidal $\sp$-Quillen adjunction.
\end{proposition}

\proof
By definition,
$U$ preserves weak equivalences and cofibrations,
and hence the pair $(U,V)$ is a Quillen adjunction.
Recall that we have the symmetric monoidal left Quillen functor
${\rm triv}:\sp\to \sp(G)_k$.
Since the composition $U\circ{\rm triv}$
is the identity functor on $\sp$,
the functor $U$ is a symmetric monoidal
left $\sp$-Quillen functor.
\qqq

\subsection{Filtered colimits in $\sp(G)_k$}

In this subsection 
we shall show that any filtered colimit preserves
weak equivalences in $\sp(G)_k$.
This follows from \cite[Prop.~4.1]{Raptis-Rosicky}
(see also \cite[Prop.~7.3]{Dugger1}
and \cite[Lem.~1.6]{Ching-Riehl}).

\begin{proposition}
\label{prop:filteredcolimit-preserves-equivalences}
Let $\Lambda$ be a filtered category,  
and let $F: X\to Y$ be a natural transformation 
of functors from $\Lambda$ to $\sp(G)_k$. 
If $F(\lambda): X(\lambda)\to Y(\lambda)$
is a $k$-local equivalence for all $\lambda\in \Lambda$,
then the induced map 
on colimits
\[ \subrel{\lambda\in\Lambda}{\rm colimit} X\longrightarrow\
   \subrel{\lambda\in\Lambda}{\rm colimit} Y \]
is also a $k$-local equivalence.
\end{proposition}

\proof
We shall apply 
\cite[Prop.~4.1]{Raptis-Rosicky}
for $\sp(G)_k$.
By Proposition~\ref{prop:spGk-combinatorial},
$\sp(G)_k$ is a combinatorial model category.
Hence we have to show that
there exists a generating set of cofibrations 
between compact objects for $\sp(G)_k$.
Let 
$I$ be the set of maps of the form
$(\partial\Delta^m\times G/N)_+\to
        (\Delta^m\times G/N)_+$ in $\sset(G)_*$,
where $m\ge 0$ and $N$ is an open subgroup of $G$.
Let ${\rm Ev}_n: \sp(G)\to \sset(G)_*$ be
the evaluation functor which assigns 
to a discrete symmetric $G$-spectrum $X$
the $n$th pointed simplicial discrete $G$-set $X_n$,
and let
$F_n: \sset(G)_*\to \sp(G)$ be its left adjoint.
We can take 
$I^{\Sigma}=\bigcup_{n\ge 0}F_n(I)$
as a set of generating cofibrations of $\sp(G)_k$. 
We can verify that 
$F_n((\partial\Delta^m\times G/N)_+)$
and $F_n((\Delta^m\times G/N)_+)$
are compact objects in $\sp(G)_k$
for all $n\ge 0, m\ge 0$ and $N$
since $\partial\Delta^m$ and $\Delta^m$ are compact objects in 
the category of simplicial sets.
This completes the proof.
\qqq

\begin{corollary}
For any functor $X$ from a filtered category 
$\Lambda$ to $\sp(G)_k$,
the canonical map
\[ \subrel{\lambda\in\Lambda}{\rm hocolim}
   X\longrightarrow\
   \subrel{\lambda\in\Lambda}{\rm colim} X\]
in ${\rm Ho}(\sp(G)_k)$
is an isomorphism.
\end{corollary}

\proof
Let $X'\to X$ be a cofibrant replacement
in the functor category ${\rm Fun}(\Lambda,\sp(G)_k)$
with the projective model structure.
This model structure exists 
by \cite[Prop.~A.2.8.2]{Lurie}
since $\sp(G)_k$ is a combinatorial model category
by Proposition~\ref{prop:spGk-combinatorial}.
The homotopy colimit ${\rm hocolim}\,X$
is represented by ${\rm colim}\,X'$.
The map ${\rm colim}\,X'\to {\rm colim}\,X$
induced on the colimits
is a $k$-local equivalence  
by Proposition~\ref{prop:filteredcolimit-preserves-equivalences}.
\qqq

\section{Embeddings of modules into $\sp(G)_k$}
\label{section:embedding_model_module}

In this section we discuss embeddings of 
certain full subcategories of module objects in $\spk$
into categories of module objects in $\sp(G)_k$.
Let $A$ be a monoid object in $\spk$.
We regard $A$ as a monoid object in $\sp(G)_k$
with trivial $G$-action.
For a map $\varphi: A\to B$ of monoid objects in $\sp(G)_k$,
we show that a certain full subcategory of 
${\rm Ho}(\mod_A(\spk))$ can be embedded into
${\rm Ho}(\mod_B(\sp(G)_k))$ as 
an ${\rm Ho}(\sp)$-enriched category.
We also discuss the relationship between
the embeddings and the homotopy fixed points functors.

\subsection{Model structure on module categories}
\label{subsection:model-structure-on-module-categories}

In this subsection
we define the model structure on 
the category of module objects
in a combinatorial symmetric monoidal
$\sp$-model category satisfying the monoid axiom.
We also study the adjunction of module categories
induced by a map of monoid objects.


Let $\mathbf{C}$ be a closed symmetric monoidal category.
First,
we recall the definitions of a closed $\mathbf{C}$-module,
a closed symmetric $\mathbf{C}$-algebra,
and an adjunction between them.
A category $\mathbf{N}$
is said to be a closed $\mathbf{C}$-module if
it is enriched, tensored, and cotensored  
over $\mathbf{C}$.
For closed $\mathbf{C}$-modules $\mathbf{N}_1$ and $\mathbf{N}_2$,
an adjoint pair of functors 
$F:\mathbf{N}_1\rightleftarrows \mathbf{N}_2:G$
is said to be an adjunction of closed $\mathbf{C}$-modules
if the left adjoint $F$ respects the tensor structures on
$\mathbf{N}_1$ and $\mathbf{N}_2$ over $\mathbf{C}$.
A category $\mathbf{A}$
is said to be a closed symmetric $\mathbf{C}$-algebra
if $\mathbf{A}$ is a closed symmetric monoidal category
equipped with a strong symmetric monoidal left adjoint functor
$i:\mathbf{C}\to\mathbf{A}$.
Let $\mathbf{A}_1$ and $\mathbf{A}_2$ be
closed symmetric $\mathbf{C}$-algebras 
equipped with strong symmetric monoidal left 
adjoints $i_1:\mathbf{C}\to\mathbf{A}_1$ 
and $i_2:\mathbf{C}\to\mathbf{A}_2$, respectively.
An adjunction 
of closed symmetric $\mathbf{C}$-algebras
between $\mathbf{A}_1$ and $\mathbf{A}_2$
is an adjoint pair of functors  
$F:\mathbf{A}_1\rightleftarrows \mathbf{A}_2:G$,
where $F$ is a strong symmetric monoidal functor,
together with a symmetric monoidal natural isomorphism
between $F\circ i_1$ and $i_2$
(see \cite[\S4.1]{Hovey} for these concepts).

In this subsection
we let $\mathbf{M}$ be a combinatorial symmetric monoidal 
$\sp$-model category with tensor product $\otimes$
satisfying the monoid axiom.
For a monoid object $R$ in $\mathbf{M}$,
we denote by $\mod_R(\mathbf{M})$ the category of 
left $R$-module objects in $\mathbf{M}$ and maps
between them.
The category $\mod_R(\mathbf{M})$ 
is a closed $\sp$-module
for a monoid object $R$ in $\mathbf{M}$.
If $R$ is a commutative monoid object,
then $\mod_R(\mathbf{M})$ is a closed symmetric 
monoidal category
with tensor product $\otimes_R$ and unit object $R$,
and furthermore,
we can regard 
$\mod_R(\mathbf{M})$ as a closed symmetric $\sp$-algebra
by 
the functor $R\otimes(-):\sp\to\mod_R(\mathbf{M})$ 
that is a strong symmetric monoidal 
left adjoint functor.

Let $R$ be a monoid object 
in $\mathbf{M}$.
A map $f: M\to N$ in $\mod_R(\mathbf{M})$
is said to be 
\begin{itemize}
\item 
a weak equivalence if it is a weak equivalence 
in $\mathbf{M}$, 
\item 
a fibration if it is a fibration in $\mathbf{M}$, and 
\item 
a cofibration if it has the left lifting property
with respect to all maps which are
both fibrations and weak equivalences.
\end{itemize}
With these definitions,
$\mod_R(\mathbf{M})$ is a model category
by \cite[Thm.~4.1]{Schwede-Shipley}.
Note that the unit object $R$ is cofibrant
in $\mod_R(\mathbf{M})$.

Let $\mathbf{C}$ be a symmetric monoidal
model category.
Now we recall the definition of a $\mathbf{C}$-model category
and a $\mathbf{C}$-Quillen adjunction
between $\mathbf{C}$-model categories
(see \cite[\S4.2]{Hovey}).
A model category $\mathbf{N}$ is said to be a $\mathbf{C}$-model category
if it is a closed $\mathbf{C}$-module and 
the action map
$\otimes : \mathbf{N}\times\mathbf{C} \to \mathbf{N}$
is a Quillen bifunctor,
that is,
for any cofibration $f: U\to V$ in $\mathbf{N}$
and any cofibration $g: W\to X$ in $\mathbf{C}$,
the induced map
\[  
   (V\otimes W)\coprod_{U\otimes W}(U\otimes X) 
   \longrightarrow V\otimes X \]
is a cofibration in $\mathbf{N}$ that is trivial
if either $f$ or $g$ is
(see \cite[Def.~4.2.1]{Hovey} for the definition
of a Quillen bifunctor).
A $\mathbf{C}$-Quillen adjunction 
between $\mathbf{C}$-model categories
is a Quillen adjunction
which is also an adjunction
of closed $\mathbf{C}$-modules. 

We can verify that
$\mod_R(\mathbf{M})$ is a $\sp$-model category
for any monoid object $R$ in $\mathbf{M}$.
If $R$ is a commutative monoid object in $\mathbf{M}$,
then $\mod_R(\mathbf{M})$ is a symmetric monoidal
model category satisfying the monoid axiom
by \cite[Thm.~4.1]{Schwede-Shipley}, and 
furthermore, since 
the strong symmetric monoidal functor 
$R\otimes(-):\sp\to \mod_R(\mathbf{M})$
is a left Quillen functor,
$\mod_R(\mathbf{M})$ is a symmetric 
monoidal $\sp$-model category
(see the paragraph just before 
Theorem~\ref{thm:weak-stable-symmetric-monoidal}
for the definition of a symmetric monoidal
$\mathbf{C}$-model category
for a symmetric monoidal model category $\mathbf{C}$).

Now we consider an adjunction between categories
of module objects induced by
a map of monoid objects in $\mathbf{M}$.
Let $\varphi:A\to B$ be a map of monoid objects in $\mathbf{M}$.
We have an adjoint pair of functors
\[ B\otimes_A (-): \mod_A(\mathbf{M})\rightleftarrows
                  \mod_B(\mathbf{M}): \varphi^*,\]
where $\varphi^*$ is the restriction of scalars functor.
We can verify that $\varphi^*$ preserves
fibrations and weak equivalences.
Hence the pair $(B\otimes_A(-),\varphi^*)$
is a Quillen adjunction.
We have natural isomorphisms
$B\otimes_A(M\otimes K)\cong (B\otimes_AM)\otimes K$
in $\mod_B(\mathbf{M})$
for $M\in\mod_A(\mathbf{M})$ and $K\in\sp$
which reduce to the canonical isomorphism
when $K$ is the sphere spectrum, 
and are 
compatible with the associativity isomorphisms
in $\mod_A(\mathbf{M})$ and $\mod_B(\mathbf{M})$
with respect to the tensor structures over $\sp$.
Hence we obtain the following lemma.

\begin{lemma}
\label{lemma:(Bwedge_A(-),varphi)-adjunction}
The adjoint pair $(B\otimes_A(-),\varphi^*)$
is a $\sp$-Quillen adjunction.
\end{lemma}

Now suppose that $A$ and $B$ are commutative monoid objects
in $\mathbf{M}$ and $\varphi: A\to B$
is a map of commutative monoid objects.
The left Quillen functor 
$B\otimes_A(-): \mod_A(\mathbf{M})\to \mod_B(\mathbf{M})$
is a strong symmetric monoidal functor
between the symmetric monoidal categories.
Since $A$ is cofibrant in $\mod_A(\mathbf{M})$,
we see that $B\otimes_A(-)$ is a symmetric monoidal
left Quillen functor.
Furthermore,
$B\otimes_A(-)$ 
is a symmetric monoidal left $\sp$-Quillen functor
since there are natural isomorphisms
$B\otimes_A(A\otimes K)\cong B\otimes K$
for $K\in\sp$.
Hence we obtain the following lemma.

\begin{lemma}
If $A$ and $B$ are commutative monoid objects
in $\mathbf{M}$ and $\varphi: A\to B$ is a map
of commutative monoid objects,
then the pair $(B\otimes_A(-),\varphi^*)$
is a symmetric monoidal $\sp$-Quillen adjunction.
\end{lemma}

\subsection{Embeddings of module categories}
\label{subsection:adjunction-model-category}

In this subsection we formulate embeddings
of module categories in combinatorial symmetric 
monoidal $\sp$-model categories satisfying the monoid axiom.
We show that a certain full subcategory
of a category of module objects  
can be embedded into another 
category of module objects 
as an enriched category
over the stable homotopy category of spectra.   

First,
we formulate the setting.
Let $\mathbf{M}$ and $\mathbf{N}$
be combinatorial symmetric monoidal 
$\sp$-model categories 
satisfying the monoid axiom.
We suppose that we have a symmetric
monoidal $\sp$-Quillen adjunction
$i: \mathbf{M}\rightleftarrows \mathbf{N}:j$.
We take monoid objects $A$ in $\mathbf{M}$
and $B$ in $\mathbf{N}$.
We suppose that there is
a morphism of monoid objects 
\[ \varphi: i(A)\longrightarrow B.\]
By Lemma~\ref{lemma:(Bwedge_A(-),varphi)-adjunction},
the morphism $\varphi$ induces a $\sp$-Quillen adjunction
$(B\otimes_{i(A)}(-),\varphi^*)$.
Since $i$ is strong symmetric monoidal,
it induces a functor
$i:\mod_A(\mathbf{M})\to \mod_{i(A)}(\mathbf{N})$.
We see that $j$ induces 
a right adjoint $j: \mod_{i(A)}(\mathbf{N})\to \mod_A(\mathbf{M})$
to $i$.
Composing these two adjunctions,
we obtain an adjoint pair of functors
\[ {\rm Ex}: \mod_A(\mathbf{M})\rightleftarrows \mod_B(\mathbf{N}) : 
   {\rm Re}.\]

\begin{lemma}
\label{lemma:module-adjunction-lemma}
The pair $({\rm Ex},{\rm Re})$ is a $\sp$-Quillen adjunction.  
\end{lemma}

\proof
By definition of the model structures,
$\varphi^*$ and $j$ are right Quillen functors.
Hence the composition ${\rm Re}=j\circ\varphi^*$
is also a right Quillen functor.
Since the pair $(i,j)$ is 
a symmetric monoidal $\sp$-Quillen adjunction,
we have natural isomorphisms
$i(M\otimes K)\cong i(M)\otimes K$
for $M\in\mathbf{M}$ and $K\in\sp$.
These isomorphisms induce
natural isomorphisms
${\rm Ex}(M\otimes K)\cong {\rm Ex}(M)\otimes K$
for $M\in\mod_A(\mathbf{M})$ and $K\in\sp$
which reduce to the canonical isomorphism when 
$K$ is the sphere spectrum,
and are compatible with the associativity 
isomorphisms in $\mod_A(\mathbf{M})$
and $\mod_B(\mathbf{N})$
with respect to the tensor structures
over $\sp$.
This completes the proof.
\qed

\bigskip

The homotopy categories $\Ho(\mod_A(\mathbf{M}))$
and $\Ho(\mod_B(\mathbf{N}))$ are
closed ${\rm Ho}(\sp)$-modules 
by \cite[Thm.~4.3.4]{Hovey}.
Furthermore,
$\Ho(\mod_A(\mathbf{M}))$
and $\Ho(\mod_B(\mathbf{N}))$ are
triangulated categories
since $\mod_A(\mathbf{M})$ and
$\mod_B(\mathbf{N})$
are stable model categories
by \cite[Lem.~3.5.2]{Schwede-Shipley2}.
By Lemma~\ref{lemma:module-adjunction-lemma},
we obtain an adjunction of 
closed $\Ho(\sp)$-modules
\[ \mathbb{L}{\rm Ex}: \Ho(\mod_A(\mathbf{M}))\rightleftarrows
   \Ho(\mod_B(\mathbf{N})):\mathbb{R}{\rm Re},\]
where $\mathbb{L}{\rm Ex}$ is the total left derived
functor of ${\rm Ex}$ and 
$\mathbb{R}{\rm Re}$ is the total right derived
functor of ${\rm Re}$.
Note that $\mathbb{L}{\rm Ex}$ and $\mathbb{R}{\rm Re}$ are
exact functors between the triangulated categories.

Now we consider the case where 
$A\in \mathbf{M}$ and $B\in\mathbf{N}$ 
are commutative monoid objects and
$\varphi: i(A)\to B$ is a map
of commutative monoid objects.
In this case $\mod_A(\mathbf{M})$
and $\mod_B(\mathbf{N})$
are symmetric monoidal $\sp$-model categories.
Since the left adjoint ${\rm Ex}$
of the Quillen adjunction
${\rm Ex}: \mod_A(\mathbf{M})\rightleftarrows
\mod_B(\mathbf{N}):{\rm Re}$ is strong symmetric monoidal,
the pair $({\rm Ex},{\rm Re})$ is a symmetric monoidal
Quillen adjunction.
Furthermore,
since there are natural isomorphisms
$B\otimes_{i(A)}i(A\otimes K)\cong
B\otimes K$ for $K\in\sp$,
we obtain the following lemma.

\begin{lemma}
\label{lemma:symmetric-monoidal-module-adjunction-lemma}
If $A\in\mathbf{M}$ and $B\in\mathbf{N}$ are 
commutative monoid objects and $\varphi:i(A)\to B$
is a map of commutative monoid objects,
then the pair $({\rm Ex},{\rm Re})$
is a symmetric monoidal $\sp$-Quillen adjunction.
\end{lemma}

If $A$ and $B$ are commutative monoid objects,
then the homotopy categories
${\rm Ho}(\mod_A(\mathbf{M}))$
and ${\rm Ho}(\mod_B(\mathbf{N}))$
are closed symmetric monoidal
${\rm Ho}(\sp)$-algebras
(see \cite[\S4.2]{Hovey}).
Furthermore,
if $\varphi: i(A)\to B$ is a map of commutative
monoid objects,
then the induced pair
$(\mathbb{L}{\rm Ex},\mathbb{R}{\rm Re})$
of functors
is an adjunction of symmetric monoidal ${\rm Ho}(\sp)$-algebras
(see the second paragraph of 
\S\ref{subsection:model-structure-on-module-categories}
for the definition of these concepts).

Now suppose $A$ and $B$ are monoid objects
and $\varphi: i(A)\to B$ is a map of monoid objects.
We shall define a full subcategory $\mathbf{T}$
of $\Ho(\mod_A(\mathbf{M}))$
and show that $\mathbf{T}$
can be embedded into $\Ho(\mod_B(\mathbf{N}))$
as an $\Ho(\sp)$-enriched category
through the functor $\mathbb{L}{\rm Ex}$.
Let $\mathbf{T}$ be the full subcategory 
of $\Ho(\mod_A(\mathbf{M}))$ consisting of $X$
such that the unit map 
$X\to \mathbb{R}{\rm Re}\mathbb{L}{\rm Ex}(X)$
is an isomorphism
\[ \mathbf{T}=\{X\in\Ho(\mod_A(\mathbf{M}))|\ 
   X\stackrel{\cong}{\longrightarrow}
   \mathbb{R}{\rm Re}\mathbb{L}{\rm Ex}(X)\}.\]
It is easy to see that 
$\mathbf{T}$ is a thick subcategory
of $\Ho(\mod_A(\mathbf{M}))$.

\begin{proposition}
\label{prop:fundamental-proposition-module-embedding}
The restriction of $\mathbb{L}{\rm Ex}$ 
to $\mathbf{T}$
is fully faithful as an $\Ho(\sp)$-enriched functor.
\end{proposition}

\proof
This follows from the natural isomorphism
\[ \mathbb{R}{\rm Map}_{\mod_B(\mathbf{N})}
   (\mathbb{L}{\rm Ex}(X),\mathbb{L}{\rm Ex}(Y))\cong
   \mathbb{R}{\rm Map}_{\mod_A(\mathbf{M})}
   (X,\mathbb{R}{\rm Re}\mathbb{L}{\rm Ex}(Y))\]
in $\Ho(\sp)$,
where $\mathbb{R}{\rm Map}_{\mod_A(\mathbf{M})}(-,-)$
and $\mathbb{R}{\rm Map}_{\mod_B(\mathbf{N})}(-,-)$
are derived mapping spaces of $\mod_A(\mathbf{M})$
and $\mod_B(\mathbf{N})$, respectively.
\qqq


\subsection{Homotopy fixed points functor}
\label{subsection:homotopy-fixed-points}

In this subsection
we shall discuss the relationship between
the homotopy fixed points functors and 
the embeddings considered 
in \S\ref{subsection:adjunction-model-category}.   
We show that the full subcategory $\mathbf{T}$
contains all dualizable objects
in some appropriate settings.

First, we recall the definition
of homotopy fixed points functors
(see \cite[\S\S3 and 5]{Behrens-Davis}).
By Proposition~\ref{prop:(i,(-)^G-adjunction-localizaed-version)},
we have a symmetric monoidal
$\sp$-Quillen adjunction
\[ {\rm triv}: \spk\rightleftarrows \sp(G)_k: (-)^G.\]
Let $A$ be a monoid object in $\spk$
and we regard $A$ as a monoid object in $\sp(G)_k$
with trivial $G$-action.
Let $\varphi: A\to B$
be a map of monoid objects 
in $\sp(G)_k$.
By Lemma~\ref{lemma:module-adjunction-lemma},
we have a $\sp$-Quillen adjunction 
\[ {\rm Ex}: \mod_A(\spk)\rightleftarrows \mod_B(\sp(G)_k): {\rm Re}.\]
This induces an adjunction 
\[ \mathbb{L}{\rm Ex}: \Ho(\mod_A(\spk))\rightleftarrows
                \Ho(\mod_B(\sp(G)_k)):\mathbb{R}{\rm Re}\]
of closed $\Ho(\sp)$-modules. 
We define a homotopy fixed points functor $(-)^{hG}$
to be the total right derived functor of ${\rm Re}$:
\[ (-)^{hG}=\mathbb{R}{\rm Re}: \Ho(\mod_B(\sp(G)_k))\longrightarrow
            \Ho(\mod_A(\spk)).\] 

Although the definition of the homotopy fixed points spectrum 
$X^{hG}$ depends on the map $\varphi$,
we shall show that the homotopy type of 
the underlying spectrum
of $X^{hG}$ is independent of $\varphi$
and also agrees with the underlying notion of 
homotopy fixed points on $\sp(G)_k$.
There is a diagram
\[ \begin{array}{ccc}
     \mod_B(\sp(G)_k)& 
     \stackrel{\rm Re}{\hbox to 15mm{\rightarrowfill}}&
     \mod_A(\spk)\\[2mm]
     \mbox{$\scriptstyle F_B$}\bigg\downarrow
     \phantom{\mbox{$\scriptstyle \pi_B$}}
     & & 
     \phantom{\mbox{$\scriptstyle \pi_A$}}\bigg\downarrow
     \mbox{$\scriptstyle F_A$}\\[2mm]
     \sp(G)_k& 
     \stackrel{(-)^G}{\hbox to 15mm{\rightarrowfill}}&
     \spk,\\
   \end{array}\]
where $F_A: \mod_A(\spk)\to \spk$ and
$F_B: \mod_B(\sp(G)_k)\to \sp(G)_k$ are 
forgetful functors.
This diagram is commutative
on the nose
since the functor ${\rm Re}: \mod_B(\sp(G)_k)\to \mod_A(\spk)$
is the composition $(-)^G\circ \varphi^*$,
where $\varphi^*:\mod_B(\sp(G)_k)\to \mod_A(\sp(G)_k)$
is the restriction of scalars functor and
$(-)^G: \mod_A(\sp(G)_k)\to \mod_A(\spk)$
is the fixed points functor on $\mod_A(\sp(G)_k)$
induced by the fixed points functor on $\sp(G)_k$.
Since $F_A$ and $F_B$ preserve weak equivalences,
they induce functors
$hF_A: \Ho(\mod_A(\spk))\to\Ho(\spk)$
and 
$hF_B:\Ho(\mod_B(\sp(G)_k))\to\Ho(\sp(G)_k)$
between the homotopy categories,
respectively.
Furthermore, since $F_A$ and $F_B$ preserve
fibrations,
we obtain a natural isomorphism
\[ hF_A(X^{hG})\stackrel{\cong}{\longrightarrow}
      (hF_BX)^{hG}\]
in $\Ho(\spk)$
for any $X\in{\rm Ho}(\mod_B(\sp(G)_k))$.
This means the homotopy type of the underlying spectrum 
of $X^{hG}$ is independent of $\varphi$
and also agrees with the underlying notion of 
homotopy fixed points on $\sp(G)_k$.

Since $A$ is cofibrant in $\mod_A(\spk)$,
we have an isomorphism $\mathbb{L}{\rm Ex}(A)\cong B$.
Hence the unit of the adjunction 
$(\mathbb{L}{\rm Ex},\mathbb{R}{\rm Re})$
gives a map $A\to B^{hG}$ in $\Ho(\mod_A(\spk))$.
Recall that $\mathbf{T}$ is the full subcategory
of $\Ho(\mod_A(\spk))$ consisting of $X$
such that the unit map 
$X\to \mathbb{R}{\rm Re}\mathbb{L}{\rm Ex}(X)$
is an isomorphism.
Hence $A\in\mathbf{T}$ if and only if
the map $A\to B^{hG}$ is an isomorphism.

Now we suppose $A\in \spk$
and $B\in\sp(G)_k$ are commutative monoid objects
and $\varphi: A\to B$ is a map of commutative monoid
objects.
In this case 
the pair
$\mathbb{L}{\rm Ex}: {\rm Ho}(\mod_A(\spk))\rightleftarrows
{\rm Ho}(\mod_B(\sp(G)_k)):\mathbb{R}{\rm Re}$
of functors
is an adjunction of symmetric monoidal ${\rm Ho}(\sp)$-algebras.
In particular,
${\rm Ho}(\mod_A(\spk))$ and ${\rm Ho}(\mod_B(\sp(G))_k)$
are symmetric monoidal categories
and $\mathbb{L}{\rm Ex}$
is a strong symmetric monoidal functor.
We show that $\mathbf{T}$ contains
all dualizable objects in ${\rm Ho}(\mod_A(\spk))$
if $A\in\mathbf{T}$.

\begin{proposition}
\label{prop:T-homotopy-fixedpoint-condition}
We assume that $A\in\spk$ and $B\in\sp(G)_k$
are commutative monoid objects and $\varphi:
A\to B$ is a map of commutative monoid objects. 
If $A\to B^{hG}$ is an isomorphism,
then $\mathbf{T}$ contains
all dualizable objects in $\Ho(\mod_A(\spk))$.
\end{proposition}

\proof
We put $\mathbf{M}=\mod_A(\spk)$ and
$\mathbf{N}=\mod_B(\sp(G)_k)$.
The homotopy categories ${\rm Ho}(\mathbf{M})$
and ${\rm Ho}(\mathbf{N})$
are symmetric monoidal categories
with tensor products $\otimes^{\mathbb{L}}_A$
and $\otimes^{\mathbb{L}}_B$, respectively. 
We let $X$ be a dualizable object in $\Ho(\mathbf{M})$
and denote by $DX$ its dual.
Since $\mathbb{L}{\rm Ex}$ is strong symmetric monoidal
and $\mathbb{L}{\rm Ex}(A)\cong B$,
$\mathbb{L}{\rm Ex}(X)$ is a dualizable object
in $\Ho(\mathbf{N})$ and
its dual is $\mathbb{L}{\rm Ex}(DX)$.
For any $W\in\Ho(\mathbf{M})$, we have 
a natural isomorphism
\[ [W,\mathbb{R}{\rm Re}\mathbb{L}{\rm Ex}(X)]_{\mathbf{M}}\cong
    [\mathbb{L}{\rm Ex}(W)\otimes^{\mathbb{L}}_B\mathbb{L}{\rm Ex}(DX),
    \mathbb{L}{\rm Ex}(A)]_{\mathbf{N}}. \]
We have 
$\mathbb{L}{\rm Ex}(W)\otimes^{\mathbb{L}}_B\mathbb{L}{\rm Ex}(DX)\cong
\mathbb{L}{\rm Ex}(W\otimes^{\mathbb{L}}_ADX)$.
This implies the following isomorphism
\[ [\mathbb{L}{\rm Ex}(W)\otimes^{\mathbb{L}}_B\mathbb{L}{\rm Ex}(DX),
    \mathbb{L}{\rm Ex}(A)]_{\mathbf{N}}\cong
   [W,\mathbb{R}{\rm Re}\mathbb{L}{\rm Ex}(A)\otimes^{\mathbb{L}}_A
    X]_{\mathbf{M}}.\]
By the Yoneda lemma,
we obtain
$\mathbb{R}{\rm Re}\mathbb{L}{\rm Ex}(X)\cong
\mathbb{R}{\rm Re}\mathbb{L}{\rm Ex}(A)\otimes^{\mathbb{L}}_AX$.
By the assumption,
we have $A\cong \mathbb{R}{\rm Re}\mathbb{L}{\rm Ex}(A)$,
and hence $\mathbb{R}{\rm Re}\mathbb{L}{\rm Ex}(A)
\otimes^{\mathbb{L}}_AX\cong X$.
This shows $X\cong \mathbb{R}{\rm Re}\mathbb{L}{\rm Ex}(X)$.
\qqq

\begin{remark}\rm
\label{remark:consistent-G-Galois-extension}
We note that $A\to B^{hG}$ is an isomorphism
under some conditions
if $B$ is a consistent $k$-local $G$-Galois extension.

First, we recall the definition of 
$k$-local Amitsur derived completion
(see, for example, \cite[Def.~8.2.1]{Rognes}).
Let $A$ be a $k$-local cofibrant commutative symmetric
ring spectrum.
For a cofibrant commutative $A$-algebra $C$,
the $k$-local Amitsur derived completion
$A_{k,C}^{\wedge}$ is the homotopy limit of the 
cosimplicial spectrum $L_kC^{\wedge_A \bullet+1}$ given by
\[ L_kC^{\wedge_A n+1}=L_k(\overbrace{C\wedge_A\cdots \wedge_A C}^{n+1}) \]
with the usual cosimplicial structure,
where 
$L_k$ is the localization functor with respect to $k$.

The $k$-local $G$-Galois extension $\varphi: A\to B$ 
in the sense of \cite[Def.~6.2.1]{Behrens-Davis} 
is said to be consistent
if the coaugmentation
of the $k$-local Amitsur derived completion
$A\to A_{k,B}^{\wedge}$ is an equivalence
(see \cite[Def.~1.0.4(1)]{Behrens-Davis}).

We suppose the localization functor
$L_k$ is given as a composite of
two localization functors $L_ML_T$,
where $L_T$ is a smashing localization and 
$L_M$ is a localization with respect to a finite spectrum $M$ 
(cf.~{\cite[Assumption~1.0.3]{Behrens-Davis}}).
Furthermore, we suppose 
$G$ has finite virtual cohomological dimension.
If $\varphi: A\to B$ is a consistent $k$-local 
$G$-Galois extension,
then $\varphi$ induces an isomorphism
$A\stackrel{\cong}{\to}B^{hG}$
by \cite[Prop.~6.1.7(3) and Cor.~6.3.2]{Behrens-Davis}.
\end{remark}

\subsection{Finite Galois extensions}

In this subsection,
as an example,
we consider the case associated to
a finite Galois extension of symmetric spectra.
Now suppose that $G$ is a finite group
and $k$ is an arbitrary symmetric spectrum.
We assume that $\varphi: A\to B$ is a faithful $k$-local 
$G$-Galois extension of symmetric spectra
(see \cite[Ch.~4]{Rognes}).
We have a symmetric monoidal $\sp$-Quillen  adjunction 
${\rm Ex}: \mod_A(\spk)\rightleftarrows
             \mod_B(\sp(G)_k): {\rm Re}$
and its derived adjunction
$\mathbb{L}{\rm Ex}: {\rm Ho}(\mod_A(\spk))
\rightleftarrows
{\rm Ho}(\mod_B(\sp(G)_k)): \mathbb{R}{\rm Re}$.
In this case
we shall show that 
the full subcategory $\mathbf{T}$
is the whole category ${\rm Ho}(\mod_A(\spk))$,
where $\mathbf{T}$ 
consists of objects $X$ in ${\rm Ho}(\mod_A(\spk))$
such that the unit map $X\to \mathbb{R}{\rm Re}\mathbb{L}{\rm Ex}(X)$
is an isomorphism.
Furthermore, we shall show that the adjoint
pair $({\rm Ex},{\rm Re})$
is a Quillen equivalence.

First, we recall that an $A$-module object $M$
is said to be $k$-locally faithful if
$M\wedge_A N\cong 0$ implies $N\cong 0$
in ${\rm Ho}(\mod_A(\spk))$.
A $k$-local $G$-Galois extension $\varphi: A\to B$
is said to be faithful if $UB$ is a $k$-locally faithful
$A$-module,
where $U: \sp(G)_k\to \spk$
is the forgetful functor.
By \cite[Prop.~6.2.1]{Rognes}, 
if $\varphi: A\to B$
is a $k$-local $G$-Galois extension,
then $UB$ is a dualizable object
in ${\rm Ho}(\mod_A(\spk))$.

For any $N\in {\rm Ho}(\mod_A(\spk))$,
we have the unit map 
$N\to \mathbb{R}{\rm Re}\mathbb{L}{\rm Ex}(N)$.
Since $\mathbb{R}{\rm Re}=(-)^{hG}$
is the homotopy fixed points functor,
we can identify this map
with the map $N\to (B\wedge_A N)^{hG}$.
We show that the map $N\to B(\wedge_A N)^{hG}$
is an isomorphism for any $N\in {\rm Ho}(\mod_A(\spk))$.
This means the full subcategory $\mathbf{T}$
is the whole category ${\rm Ho}(\mod_A(\spk))$.

\begin{lemma}
\label{lem:G-finite-T-whole}
Suppose that $G$ is a finite group and $\varphi: A\to B$
is a faithful $k$-local $G$-Galois extension.
The map $N\to (B\wedge_A N)^{hG}$
is an isomorphism for any $N\in {\rm Ho}(\mod_A(\spk))$,
and hence 
$\mathbf{T}={\rm Ho}(\mod_A(\spk))$.
\end{lemma}

\proof
Since $UB$ is $k$-locally faithful over $A$,
it suffices to show that
the map $UB\wedge_A N\to UB\wedge_A (B\wedge_A N)^{hG}$
is an isomorphism in ${\rm Ho}(\mod_A(\spk))$.
Since $UB$ is a dualizable object, 
we have an isomorphism
$UB\wedge_A (B\wedge_A N)^{hG}\cong
(UB\wedge_A B\wedge_A N)^{hG}$
in ${\rm Ho}(\mod_A(\spk))$.
By the assumption that $\varphi: A\to B$
is a $k$-local $G$-Galois extension,
we have an equivalence
$UB\wedge_AB\simeq {\rm Map}(G,UB)$.
This induces an equivalence
$UB\wedge_AB\wedge_AN\simeq {\rm Map}(G,UB\wedge_A,N)$.
Hence we see that
$(UB\wedge_A B\wedge_A N)^{hG}\simeq UB\wedge_A N$
and the map $UB\wedge_A N\to UB\wedge_A (B\wedge_A N)^{hG}$
is an isomorphism in ${\rm Ho}(\mod_A(\spk))$.
This completes the proof.
\qqq

Next, we would like to show that
the adjunction $({\rm Ex},{\rm Re})$
is actually a Quillen equivalence.
In order to show that
the following lemma is useful.

\begin{lemma}
\label{lem:Z-hg-trivial-implies-Z-trivial}
Suppose that $G$ is a finite group and $\varphi: A\to B$
is a faithful $k$-local $G$-Galois extension.
For any $Z\in{\rm Ho}(\mod_B(\sp(G)_k))$,
if $Z^{hG}$ is trivial,
then $Z$ is trivial.
\end{lemma}

\proof
We may assume that $Z$ is an object 
in $\mod_B(\sp(G)_k)$ that is both fibrant and cofibrant.
We have a map
$UB\wedge_A Z\to {\rm Map}(G,UZ)$ in $\sp(G)_k$,
which is the adjoint of 
the $UB$-action map $U(UB\wedge_A Z)\cong UB\wedge_A UZ
\to UZ$ of $UZ$.
We shall show that this map is a weak equivalence.
For this purpose,
it suffice to show that
$UB\wedge_A UZ\to U{\rm Map}(G,UZ)$
is a weak equivalence.
Using the fact
that $UB\wedge_AUB\simeq U{\rm Map}(G,UB)$,
we see that
$UB\wedge_A UZ\simeq U{\rm Map}(G,UZ)$
and the map $UB\wedge_A UZ\to U{\rm Map}(G,UZ)$
is a weak equivalence.

The weak equivalence
$UB\wedge_A Z\to {\rm Map}(G,UZ)$
implies that
$(UB\wedge_AZ)^{hG}\cong {\rm Map}(G,UZ)^{hG}$
in ${\rm Ho}(\mod_A(\spk))$.
Since $UB$ is a dualizable object,
we have an isomorphism 
$(UB\wedge_AZ)^{hG}\cong UB\wedge_A Z^{hG}$
in ${\rm Ho}(\mod_A(\spk))$.
Since ${\rm Map}(G,UZ)^{hG}\cong UZ$,
we see that, if $Z^{hG}$ is trivial,
$UZ$ is trivial and hence $Z$ is trivial.
\qqq

Using this lemma,
we obtain the following proposition
which says the adjunction
$({\rm Ex},{\rm Re})$ is a Quillen equivalence
if it is associated to a faithful $k$-local finite Galois extension.

\begin{proposition}
If $G$ is a finite group and
$\varphi: A\to B$ is a faithful $k$-local $G$-Galois extension,
then the adjoint pair 
\[ {\rm Ex}:\mod_A(\spk)\rightleftarrows 
\mod_B(\sp(G)_k):{\rm Re} \]
is a Quillen equivalence.
\end{proposition}

\proof
Let $X$ be a cofibrant object in $\mod_A(\spk)$
and $Y$ a fibrant object in $\mod_B(\sp(G)_k)$.
Suppose that $f: B\wedge_A X\to Y$
is a weak equivalence in $\mod_B(\sp(G)_k)$.
Then we can regard $Y$ as a fibrant replacement
of $B\wedge_A X$ and
hence $Y^G$ represents $(B\wedge_A X)^{hG}$
in ${\rm Ho}(\mod_A(\spk))$.
By Lemma~\ref{lem:G-finite-T-whole},
the adjoint map
$g:X\to Y^G$ to $f$ induces an isomorphism
$X\stackrel{\cong}{\to}(B\wedge_A X)^{hG}$ in 
the homotopy category ${\rm Ho}(\mod_A(\spk))$
and hence the map $g:X\to Y^G$ is a weak equivalence.

Conversely,
suppose that $g: X\to Y^G$ 
is a weak equivalence
in $\mod_A(\spk)$.
We have to show that
the adjoint map
$f: B\wedge_A X\to Y$ to $g$ is a weak equivalence.
We let $Z$ be the cofiber 
of the map $f$
in ${\rm Ho}(\mod_B(\sp(G)_k))$.
Note that $Y^{hG}$ is represented by $Y^G$
and $(B\wedge_A X)^{hG}\cong X$ 
by Lemma~\ref{lem:G-finite-T-whole}.
By the assumption that
$g:X\to Y^G$ is a weak equivalence,
we see that the induced map
$f^{hG}:(B\wedge_AX)^{hG}\to Y^{hG}$ is an isomorphism
in ${\rm Ho}(\mod_A(\spk))$,
and hence $Z^{hG}$ is trivial.
By Lemma~\ref{lem:Z-hg-trivial-implies-Z-trivial},
$Z$ is trivial and hence
$f$ is an isomorphism
in ${\rm Ho}(\mod_B(\sp(G)_k))$.
This shows that $f: B\wedge_AX\to Y$
is a weak equivalence.
This completes the proof.
\qqq


\subsection{The $K(n)$-local category}
\label{subsection:K(n)-local-category}

In this subsection 
we shall apply the results 
in \S\S\ref{subsection:adjunction-model-category} 
and \ref{subsection:homotopy-fixed-points}
to the $K(n)$-local category.
Let $E_n$ be the $n$th Morava $E$-theory spectrum 
and $K(n)$ the $n$th Morava $K$-theory spectrum
at a prime $p$.
We denote by $\mathbb{G}_n$ 
the extended Morava stabilizer group.
In this subsection
we shall show that the $K(n)$-local category
can be embedded as an enriched category
over the stable homotopy category of spectra
into the homotopy category of module
objects over a discrete model $F_n$
of $E_n$ in the category of $K(n)$-local
discrete symmetric $\mathbb{G}_n$-spectra.

In \cite[Def.~2.3]{Davis} 
Davis constructed a discrete $\mathbb{G}_n$-spectrum $F_n$
by using results in \cite{Devinatz-Hopkins}.
We shall review the construction of $F_n$.
Devinatz-Hopkins \cite{Devinatz-Hopkins}
gave a functorial construction
of commutative $S$-algebras $E_n^{dhU}$ 
for open subgroups $U$ of $\mathbb{G}_n$
which satisfy the desired properties of
the homotopy fixed points spectra.
The spectrum $F_n$ is defined by
\[ F_n=\ \subrel{U}{\rm colim} E_n^{dhU}, \] 
where the colimit is taken over the open subgroups
$U$ of $\mathbb{G}_n$.  
Note that $F_n$ is $E(n)$-local,
where $E(n)$ is the $n$th Johnson-Wilson spectrum at $p$,
since the Bousfield localization $L_{E(n)}$ with
respect to $E(n)$ is smashing
and $F_n$ is the colimit of the $E(n)$-local spectra $E_n^{dhU}$,
and that 
the spectrum $E_n$ is recovered by the equivalence
\[ E_n\simeq L_{K(n)}F_n\]
in the stable homotopy category of spectra
by \cite[Thm~6.3]{Davis}, 
where $L_{K(n)}$ is the Bousfield localization
with respect to $K(n)$.

In \cite[\S8.1]{Behrens-Davis} 
Behrens-Davis showed that $F_n$ can be taken as a commutative monoid 
object in $\sp(\mathbb{G}_n)$.
Furthermore,
they showed that 
$F_n$ is a consistent
$K(n)$-local $\mathbb{G}_n$-Galois extension
of $L_{K(n)}S$,
where $S$ is the sphere spectrum.
Note the $K(n)$-localization functor $L_{K(n)}$
is the composite $L_{F(n)}L_{E(n)}$ of two localization functors,
where $L_{E(n)}$ is smashing
and
$L_{F(n)}$ is the localization functor with respect
to any finite spectrum $F(n)$ of type $n$ at $p$.
By Remark~\ref{remark:consistent-G-Galois-extension},
this implies that the unit map 
$\varphi: S\to F_n$ induces an 
isomorphism $S\stackrel{\cong}{\to} (F_n)^{h\mathbb{G}_n}$ in
$\Ho(\sp_{K(n)})$
since $\mathbb{G}_n$ has finite 
virtual cohomological dimension.

We consider the adjunction
\[ \mathbb{L}{\rm Ex}: \Ho(\sp_{K(n)})\rightleftarrows
   \Ho(\mod_{F_n}(\sp(\mathbb{G}_n)_{K(n)})):\mathbb{R}{\rm Re}.\]
By \S\ref{subsection:homotopy-fixed-points},
the right adjoint $\mathbb{R}{\rm Re}$
is interpreted as 
the homotopy fixed points functor $(-)^{h\mathbb{G}_n}$.
Hence the unit map of the adjunction
is regarded as a map
$X\to ((F_n\wedge X)_{K(n)})^{h\mathbb{G}_n}$
for $X\in {\rm Ho}(\sp_{K(n)})$,
where $(-)_{K(n)}$ is 
the $K(n)$-localization functor 
on ${\rm Ho}(\mod_{F_n}(\sp(\mathbb{G}_n))$.

Now we shall show that 
the homotopy fixed points functor $(-)^{h\mathbb{G}_n}$
is compatible with the $K(n)$-localization.
 We let $M$
be an object in $\Ho(\mod_{F_n}(\sp(\mathbb{G}_n)))$
and denote by $M_{K(n)}$
the $K(n)$-localization of $M$,
which is represented by 
a fibrant replacement in $\mod_{F_n}(\sp(\mathbb{G}_n)_{K(n)})$
of a representative of $M$ in $\mod_{F_n}(\sp(\mathbb{G}_n))$.
Since $(M_{K(n)})^{h\mathbb{G}_n}$ is $K(n)$-local,
we have a natural map
\[ L_{K(n)}(M^{h\mathbb{G}_n})\longrightarrow 
    (M_{K(n)})^{h\mathbb{G}_n}\]
in the stable homotopy category of spectra.

\begin{lemma}
\label{lem:fixed-point-K(n)-localization-exchange}
The natural map 
$L_{K(n)}(M^{h\mathbb{G}_n})\to
   (M_{K(n)})^{h\mathbb{G}_n}$
is an isomorphism in the stable homotopy category of spectra
for any object $M$ in $\Ho(\mod_{F(n)}(\sp(\mathbb{G}_n)))$.
\end{lemma}

\proof
Let $F(n)$ be a finite spectrum of type $n$ at $p$.
To show that the map $L_{K(n)}(M^{h\mathbb{G}_n})\to (M_{K(n)})^{h\mathbb{G}_n}$
is an isomorphism,
it suffices to show that
$F(n)\wedge L_{K(n)}(M^{h\mathbb{G}_n})\to
F(n)\wedge (M_{K(n)})^{h\mathbb{G}_n}$
is an isomorphism
since $L_{K(n)}(M^{h\mathbb{G}_n})$
and $(M_{K(n)})^{h\mathbb{G}_n}$
are both $K(n)$-local.

First, note that $M$ is $E(n)$-local since $F_n$ is $E(n)$-local 
and $L_{E(n)}$ is smashing.
Since $M$ and $M_{K(n)}$
are $E(n)$-local,
the localization map $M\to M_{K(n)}$ 
induces an isomorphism 
$F(n)\wedge M\stackrel{\cong}{\to}
F(n)\wedge M_{K(n)}$
in ${\rm Ho}(\mod_{F_n}(\sp(\mathbb{G}_n))$.
Hence a fibrant representative of $F(n)\wedge M$
in $\mod_{F_n}(\sp(\mathbb{G}_n))$
is also fibrant in $\mod_{F_n}(\sp(\mathbb{G}_n)_{K(n)})$.
Therefore,
we see that the homotopy fixed points spectrum
$(F(n)\wedge M)^{h\mathbb{G}_n}$
taken on ${\rm Ho}(\mod_{F_n}(\sp(\mathbb{G}_n)))$
is isomorphic in the stable homotopy category of spectra
to the homotopy fixed points spectrum
$(F(n)\wedge M_{K(n)})^{h\mathbb{G}_n}$
taken on ${\rm Ho}(\mod_{F_n}(\sp(\mathbb{G}_n))_{K(n)})$.

Since $M$ is $E(n)$-local,
$M^{h\mathbb{G}_n}$ is also $E(n)$-local.
Using the facts that $F(n)$ is finite, $M^{h\mathbb{G}_n}$ is $E(n)$-local,
and $L_{K(n)}\simeq L_{F(n)}L_{E(n)}$,
we have a natural isomorphism
$F(n)\wedge L_{K(n)}(M^{h\mathbb{G}_n})\cong
(F(n)\wedge M)^{h\mathbb{G}_n}$.
On the other hand,
since $F(n)$ is finite,
we have a natural isomorphism
$F(n)\wedge (M_{K(n)})^{h\mathbb{G}_n}\cong
(F(n)\wedge M_{K(n)})^{h\mathbb{G}_n}$. 
Therefore,
we see that the map
$F(n)\wedge L_{K(n)}(M^{h\mathbb{G}_n})\to
F(n)\wedge (M_{K(n)})^{h\mathbb{G}_n}$
is an isomorphism.
This completes the proof.
\qqq

By \cite[Thm.~9.7]{Davis} and \cite[Thm.~1.1]{Davis-Torii},
the natural map
$X\to L_{K(n)}((F_n\wedge X)^{h\mathbb{G}_n})$
is an isomorphism 
for any $X\in {\rm Ho}(\sp_{K(n)})$.
Using Lemma~\ref{lem:fixed-point-K(n)-localization-exchange},
we see that 
the unit map $X\to ((F_n\wedge X)_{K(n)})^{h\mathbb{G}_n}$
is an isomorphism for any $X\in \Ho(\sp_{K(n)})$.
Hence we obtain the following theorem 
by Proposition~\ref{prop:fundamental-proposition-module-embedding}.

\begin{theorem}
\label{theorem:K(n)-local-embedding}
The functor
\[ \mathbb{L}{\rm Ex}: \Ho(\sp_{K(n)})\longrightarrow
                \Ho(\mod_{F_n}(\sp(\mathbb{G}_n)_{K(n)}))\]
is fully faithful as an $\Ho(\sp)$-enriched functor.
\end{theorem}

The theory of localizations in enriched categories
was developed by Wolff~\cite{Wolff}.
Theorem~\ref{theorem:K(n)-local-embedding} implies that 
the $K(n)$-local category
$\Ho(\sp_{K(n)})$ is an $\Ho(\sp)$-enriched coreflective subcategory of
$\Ho(\mod_{F_n}(\sp(\mathbb{G}_n))_{K(n)})$.
By \cite[Thm.~1.6]{Wolff},
we obtain the following corollary.

\begin{corollary}
Let $\mathbf{W}$
be the class of morphisms $f$
in $\Ho(\mod_{F_n}(\sp(\mathbb{G}_n)_{K(n)}))$
such that 
the induced morphism
$f^{h\mathbb{G}_n}=\mathbb{R}{\rm Re}(f)$
on homotopy fixed points spectra is an isomorphism
in $\Ho(\sp_{K(n)})$.
The $K(n)$-local category $\Ho(\sp_{K(n)})$
is equivalent to the localization of
the homotopy category
$\Ho(\mod_{F_n}(\sp(\mathbb{G}_n)_{K(n)}))$ 
with respect to $\mathbf{W}$
as an $\Ho(\sp)$-enriched category
\[ \Ho(\sp_{K(n)})\simeq \Ho(\mod_{F_n}(\sp(\mathbb{G}_n)_{K(n)}))
   [\mathbf{W}^{-1}].\] 
\end{corollary}

\begin{remark}\rm
In the forthcoming paper~\cite{torii}
we will show that 
the symmetric monoidal $\sp$-Quillen adjunction
${\rm Ex}: \sp_{K(n)}\rightleftarrows
\mod_{F_n}(\sp(\mathbb{G}_n)_{K(n)}):{\rm Re}$ 
is a Quillen equivalence.
\end{remark}


\section{Embeddings of quasi-categories 
of modules over $\ispk$}
\label{section:embedding-quasi-categories}

Let $\ispk$ be the underlying quasi-category
of the simplicial model category $\spk$.
Let $\psi:A\to E$ be a map of algebra objects in $\ispk$.
We have an adjunction of underlying quasi-categories
\[ E\wedge_A(-):\mod_A(\ispk)\rightleftarrows
   \mod_E(\ispk):\psi^*.\] 
In this section we discuss an embedding of 
certain full subcategory of $\mod_A(\ispk)$
into the quasi-category of comodules
associated to the adjunction.

\subsection{Quasi-categories of comodules}
\label{subsection:quasi-category-comodule}

In this subsection
we shall introduce a quasi-category of comodules
associated to an adjunction of quasi-categories.

First, we fix notation.
For a monoidal quasi-category $\mathcal{M}$,
we denote by ${\rm Alg}(\mathcal{M})$
the quasi-category of algebra objects
in $\mathcal{M}$.
For a quasi-category $\mathcal{N}$ 
left-tensored over $\mathcal{M}$
and an algebra object $T\in {\rm Alg}(\mathcal{M})$,
we denote by ${\rm  Mod}_{T}(\mathcal{N})$
the quasi-category of left $T$-modules in $\mathcal{N}$
(see \cite[4.1.1 and 4.1.2]{Lurie2} for these concepts).

Let $\mathcal{C}$ and $\mathcal{D}$ be quasi-categories.
We denote by ${\rm Fun}(\mathcal{C},\mathcal{D})$
the quasi-category of functors from $\mathcal{C}$
to $\mathcal{D}$
(see \cite[1.2.7]{Lurie}).
We can regard the quasi-category
${\rm End}(\mathcal{C})={\rm Fun}(\mathcal{C},\mathcal{C})$
as a monoidal quasi-category by the composition
of functors
(see \cite[4.7]{Lurie2}).
A monad on $\mathcal{C}$
is defined to be 
an algebra object of ${\rm End}(\mathcal{C})$.  
If $T$ is a monad on $\mathcal{C}$,
we can consider the quasi-category of left $T$-modules 
${\rm Mod}_T(\mathcal{C})$
in $\mathcal{C}$.

For quasi-categories $\mathcal{C}$ and $\mathcal{D}$,
the quasi-category ${\rm Fun}(\mathcal{D},\mathcal{C})$
carries a left action of the monoidal quasi-category
${\rm End}(\mathcal{C})$
by composition of functors,
and we can regard ${\rm Fun}(\mathcal{D},\mathcal{C})$
as left-tensored over ${\rm End}(\mathcal{C})$.
Thus, we can consider 
the quasi-category of left $T$-modules
${\rm Mod}_T({\rm Fun}(\mathcal{D},\mathcal{C}))$
in ${\rm Fun}(\mathcal{D},\mathcal{C})$
for a monad $T\in {\rm Alg}({\rm End}(\mathcal{C}))$.
Let $R: \mathcal{D}\to\mathcal{C}$
be a functor of quasi-categories.
An endomorphism monad of $R$
consists of a monad $T\in {\rm Alg}({\rm End}(\mathcal{C}))$ 
together with a left $T$-module 
$\overline{R}\in {\rm Mod}_T({\rm Fun}(\mathcal{D},\mathcal{C}))$
whose image in ${\rm Fun}(\mathcal{D},\mathcal{C})$
coincides with $R$,
such that the action map $a:T R\to R$
induces a weak equivalence of mapping spaces
\[ {\rm Map}_{{\rm Fun}(\mathcal{D},\mathcal{C})}(F,R)
   \stackrel{T}{\longrightarrow}
   {\rm Map}_{{\rm Fun}(\mathcal{D},\mathcal{C})}(T F,T R)
   \stackrel{a_*}{\longrightarrow} 
   {\rm Map}_{{\rm Fun}(\mathcal{D},\mathcal{C})}(T F,R)\]
for any $F\in {\rm Fun}(\mathcal{D},\mathcal{C})$
(see \cite[4.7.4]{Lurie2}). 

Now we recall the definition of an adjunction between
quasi-categories
(see \cite[Def.~5.2.2.1]{Lurie}).
Let $\mathcal{C}$ and $\mathcal{D}$
be quasi-categories.
An adjunction between $\mathcal{C}$ and $\mathcal{D}$
is a map $q:\mathcal{M}\to \Delta^1$ 
of simplicial sets which 
is both a Cartesian fibration and a coCartesian
fibration together with equivalences
$\mathcal{C}\to \mathcal{M}_{\{0\}}$
and $\mathcal{D}\to \mathcal{M}_{\{1\}}$,
where $\mathcal{M}_{\{0\}}$ and $\mathcal{M}_{\{1\}}$
are the fibers of $q$ at $\{0\}\in\Delta^1$ and $\{1\}\in \Delta^1$,
respectively
(see \cite[2.4.2]{Lurie} for the definitions
of a Cartesian fibration and a coCartesian fibration). 
In this case we let $L: \mathcal{C}\to\mathcal{D}$ 
and $R:\mathcal{D}\to \mathcal{C}$ be functors
associated to $\mathcal{M}$,
and say that $L$ is left adjoint to $R$ and 
$R$ is right adjoint to $L$.

We have a characterization of adjoint functors
of quasi-categories in terms of mapping spaces
as in classical category theory.
Suppose we have a pair of functors of quasi-categories
\[ L: \mathcal{C}\rightleftarrows \mathcal{D}: R.\]
The functor $L$ is left adjoint to $R$ if and only if
there exists a morphism
$u: {\rm id}_{\mathcal{C}}\to R L$ 
in ${\rm Fun}(\mathcal{C},\mathcal{C})$
such that the composition 
\[ {\rm Map}_{\mathcal{D}}(L(C),D)
   \stackrel{R}{\longrightarrow}
   {\rm Map}_{\mathcal{C}}(RL(C),R(D))
   \stackrel{u_C^*}{\longrightarrow}
   {\rm Map}_{\mathcal{C}}(C,R(D)) \]  
is a weak equivalence 
for any $C\in\mathcal{C}$ and $D\in\mathcal{D}$
(see \cite[Prop.~5.2.2.8]{Lurie}).

We shall introduce a quasi-category
of comodules associated to an adjunction
of quasi-categories.
For a quasi-category $\mathcal{X}$,
we denote by $\mathcal{X}^{\rm op}$
the opposite quasi-category of $\mathcal{X}$
(see \cite[1.2.1]{Lurie}).
Suppose we have an adjunction of quasi-categories
$L: \mathcal{C}\rightleftarrows \mathcal{D}: R$.
This induces an adjunction of opposite quasi-categories
\[ R^{\rm op}: \mathcal{D}^{\rm op}\rightleftarrows
   \mathcal{C}^{\rm op}: L^{\rm op}.\]
We obtain an object $L^{\rm op} R^{\rm op}$
in ${\rm End}(\mathcal{D}^{\rm op})$, and 
we can lift $L^{\rm op} R^{\rm op}$
to an endomorphism monad 
of $L^{\rm op}$ by \cite[4.7.4]{Lurie2}.
In particular,
we have a monad 
$\Theta\in{\rm Alg}({\rm End}(\mathcal{D}^{\rm op}))$
that is a lifting of $L^{\rm op}R^{\rm op}$, 
and a left $\Theta$-module
object $\overline{L^{\rm op}}\in
{\rm Mod}_{\Theta}({\rm Fun}(\mathcal{C}^{\rm op},
\mathcal{D}^{\rm op}))$ that is a lifting of $L^{\rm op}$.
We regard $\Theta$ as a comonad on $\mathcal{D}$ and
we define the quasi-category of left $\Theta$-comodules
\[ \comod_{\Theta}(\mathcal{D}) \]
to be ${\rm Mod}_{\Theta}(\mathcal{D}^{\rm op})^{\rm op}$.
Using an equivalence
${\rm Mod}_{\Theta}({\rm Fun}(\mathcal{C}^{\rm op},
\mathcal{D}^{\rm op}))\stackrel{\simeq}{\to}
{\rm Fun}(\mathcal{C}^{\rm op},
{\rm Mod}_{\Theta}(\mathcal{D}^{\rm op}))$,
we obtain an object in 
${\rm Fun}(\mathcal{C}^{\rm op},
{\rm Mod}_{\Theta}(\mathcal{D}^{\rm op}))$
corresponding to $\overline{L^{\rm op}}\in
{\rm Mod}_{\Theta}(\mathcal{C}^{\rm op},\mathcal{D}^{\rm op})$.
Hence we see that the functor $L: \mathcal{C}\to\mathcal{D}$
factors through a functor 
\[ \widetilde{L}: \mathcal{C}\longrightarrow
       \comod_{\Theta}(\mathcal{D}) \]
so that $U\widetilde{L}\simeq L$, 
where $U: \comod_{\Theta}(\mathcal{D})\to
\mathcal{D}$ is the forgetful functor. 
We say that $L$ exhibits $\mathcal{C}$ as comonadic over $\mathcal{D}$
if $\widetilde{L}$ is an equivalence of quasi-categories.

\subsection{Embeddings into
quasi-categories of comodules}

Recall that $\ispk$ is the underlying 
quasi-category of the simplicial model
category $\spk$.
Since $\spk$ is a simplicial symmetric monoidal
model category,
$\ispk$ is a symmetric monoidal quasi-category
(see \cite[4.1.3]{Lurie2}).
In this subsection 
we shall formulate embeddings
of quasi-categories of modules
into quasi-categories of comodules
in $\ispk$.

For an algebra object $A$ in $\ispk$,
we have the quasi-category of left $A$-modules 
$\mod_A(\ispk)$ in $\ispk$. 
For a map $\psi: A\to E$ of algebra objects
in $\ispk$,
we have an adjunction of quasi-categories
\[ L: \mod_A(\ispk)\rightleftarrows
   \mod_E(\ispk):R, \]
where $L=E\wedge_A(-)$ and $R=\psi^*$.
Hence we obtain a comonad
$\Theta$ on $\mod_E(\ispk)$ and
a quasi-category of left $\Theta$-comodules
\[ \comod_{(E,\Theta)}(\ispk)=\comod_{\Theta}(\mod_E(\ispk)).\]
The functor $L$ factors through a functor
\[ {\rm Coex}: 
   \mod_A(\ispk)\longrightarrow \comod_{(E,\Theta)}(\ispk).\]
so that $U{\rm Coex}\simeq L$,
where $U$ is the forgetful functor 
$U: \comod_{(E,\Theta)}(\ispk)\to \mod_E(\ispk)$.
We set
\[ \Theta'={\rm Coex}\,R.\]
The functor $\Theta'$ is
informally given by $\Theta'(X)=\Theta (X)$
with the obvious $\Theta$-comodule structure
for $X\in\mod_E(\ispk)$.
Note that $\Theta'$ is a right adjoint to 
the forgetful functor $U$.

Now we introduce a functor
$P: {\rm Comod}_{(E,\Theta)}(\ispk)\to\mod_A(\ispk)$
which is a derived functor of taking primitive elements.
The functor $P$ is related to
the derived completion defined by Carlsson in \cite{Carlsson} 
and the nilpotent completion considered by
Bousfield in \cite{Bousfield}.
For $X\in\comod_{(E,\Theta)}(\ispk)$,
we have a cosimplicial object 
\[ C^{\bullet}(R,\Theta,UX)\]
in $\mod_A(\ispk)$ by the cobar construction.
We define a functor 
\[ P:\comod_{(E,\Theta)}(\ispk)\to \mod_A(\ispk) \]
by $PX=\lim C^{\bullet}(R,\Theta,UX)$.

For $Y\in\mod_A(\ispk)$,
we have a coaugmented cosimplicial object
\[ Y\to E^{\bullet+1}Y \] 
in $\mod_A(\ispk)$
given by 
\[ E^{k+1}Y=\overbrace{E\wedge_A\cdots\wedge_AE}^{k+1}\wedge_AY\]
with the usual cosimplicial structure,
which is sometimes called the Amitsur complex
\cite{Amitsur}.
There is an equivalence of cosimplicial
objects
\[ C^{\bullet}(R,\Theta,U{\rm Coex}(Y))\simeq
   E^{\bullet+1}Y.\]

Note that the map $Y\to \lim E^{\bullet+1}Y$ 
is an analogue of the derived completion of $Y$ 
at the $A$-algebra $E$
in the sense of \cite{Carlsson}
(see Remark~\ref{remark:consistent-G-Galois-extension}
for the $k$-local Amitsur derived completion).

Furthermore, if $k=A=S$,
then the map $Y\to \lim E^{\bullet+1}Y$ 
is the $E$-nilpotent completion of $Y$ in ${\rm Ho}(\isp)$
in the sense of \cite{Bousfield},
where ${\rm Ho}(\isp)$ is the stable homotopy category
of spectra.

We shall recall the $B$-nilpotent completion of spectra
for a ring spectrum $B$ in ${\rm Ho}(\isp)$.
A spectrum $W\in {\rm Ho}(\isp)$
is said to be $B$-nilpotent
if $W$ lies in the thick ideal of ${\rm Ho}(\isp)$
generated by $B$.  
A $B$-nilpotent resolution of a spectrum $Z$
is a tower $\{W_s\}_{s\ge 0}$ under $Z$ in ${\rm Ho}(\isp)$
such that $W_s$ is $B$-nilpotent for all $s\ge 0$
and the map 
${\rm colim}_s\, {\rm Hom}_{\rm Ho(\isp)}(W_s,N)\to  
{\rm Hom}_{\rm Ho(\isp)}(Z,N)$ is an isomorphism for any
$B$-nilpotent spectrum $N$.
The $B$-nilpotent completion of $Z$
is defined to be the homotopy inverse limit
${\rm holim}_s W_s$
for any $B$-nilpotent resolution $\{W_s\}_{s\ge 0}$ of $Z$.

We can consider the Tot tower $\{{\rm Tor}^sE^{\bullet+1}Y\}_{s\ge 0}$
associated to the cosimplicial object
$E^{\bullet+1}Y$.
If $k=A=S$,
then the cofiber of each map
${\rm Tot}^{s+1}E^{\bullet+1}Y\to {\rm Tot}^{s}E^{\bullet+1}Y$
is an $E$-module,
and hence $E$-nilpotent in ${\rm Ho}(\isp)$
for all $s\ge 0$.
By induction on $s$ and the fact that 
${\rm Tot}^0E^{\bullet+1}Y=E\wedge Y$,
we see that ${\rm Tot}^sE^{\bullet+1}Y$
is $E$-nilpotent for all $s\ge 0$.
Furthermore,
we let $F^s$ be the fiber of
the map $Y\to {\rm Tot}^{s}E^{\bullet+1}Y$
for $s\ge 0$.
The induced map
$F^{s+1}\to F^s$
is null in ${\rm Ho}(\isp)$ after tensoring with $E$
for all $s\ge 0$.
This implies an isomorphism
${\rm colim}_s\, {\rm Hom}_{{\rm Ho}(\isp)}({\rm Tot}^sE^{\bullet+1}Y,N)\to
{\rm Hom}_{{\rm Ho}(\isp)}(Y,N)$ for any $E$-nilpotent spectrum $N$.
Hence ${\rm lim}\, E^{\bullet+1}Y$
is the $E$-nilpotent completion of $Y$ in ${\rm Ho}(\isp)$. 

We would like to show that
a certain full subcategory of
$\mod_A(\ispk)$ can be embed into 
$\comod_{(E,\Theta)}(\spk)$
through ${\rm Coex}$.
For this purpose,
we show that ${\rm Coex}$ has a right adjoint. 

\begin{proposition}
The functor $P$ is a right adjoint to ${\rm Coex}$ so
that we have an adjunction of quasi-categories
\[ {\rm Coex}: \mod_A(\ispk)\rightleftarrows
       \comod_{(E,\Theta)}(\ispk): P.\]
\end{proposition}

\proof
Let $\mathcal{X}$ be the full subcategory of
${\rm Comod}_{(E,\Theta)}(\ispk)$ consisting of
$X$ such that the functor
\[ {\rm Map}_{\comod_{(E,\Theta)}(\ispk)}
   ({\rm Coex}(-),X):
    \mod_A(\ispk)\to\mathcal{S} \]
is representable,
where $\mathcal{S}$ is the quasi-category
of spaces.
We denote by $\widetilde{P}(X)$
the representing object in $\mod_A(\ispk)$
for $X\in\mathcal{X}$.
In this case, $\widetilde{P}(X)$ is well-defined up to
canonical equivalence and we obtain a functor
\[ \widetilde{P}:\mathcal{X}\to\mod_A(\ispk).\]

In order to prove the proposition,
we have to show that $\mathcal{X}$
is actually the whole quasi-category
${\rm Comod}_{(E,\Theta)}(\ispk)$.
First,
we shall show that $\Theta'Z\in\mathcal{X}$
for any $Z\in \mod_E(\ispk)$.
Since $\Theta'$ is a right adjoint to $U$
and $U{\rm Coex}\simeq L$ is a left adjoint to $R$,
we see that there is a natural equivalence
\[ {\rm Map}_{\comod_{(E,\Theta)}(\ispk)}
   ({\rm Coex}(Y),\Theta'Z)\simeq
   {\rm Map}_{\mod_A(\ispk)}(Y,RZ)\]
for any $Y\in\mod_E(\ispk)$, and hence
$\Theta'Z\in \mathcal{X}$.
We note that $\widetilde{P}\Theta'Z\simeq RZ$.
   
Next, we shall show that $X\in\mathcal{X}$
for any $X\in \comod_{(E,\Theta)}(\ispk)$.
Let
$C^{\bullet}(\Theta',\Theta,UX)$
be a cosimplicial object in $\comod_{(E,\Theta)}(\ispk)$
given by the cobar construction.
By \cite[4.7.4]{Lurie2},
we see that 
$\lim C^{\bullet}(\Theta',\Theta,UX)\simeq X$.
Since $C^n(\Theta',\Theta,UX)\simeq \Theta'\Theta^nUX
\in\mathcal{X}$ for any $n$, 
we can construct a cosimplicial object
$\widetilde{P}C^{\bullet}(\Theta',\Theta,UX)$
in $\mod_A(\ispk)$. 
We can verify that 
$\lim \widetilde{P}C^{\bullet}(\Theta',\Theta,UX)$
represents the functor
${\rm Map}_{\comod_{(E,\Theta)}(\ispk)}({\rm Coex}(-),X)$,
and hence $X\in \mathcal{X}$.

Therefore, we obtain $\mathcal{X}=\comod_{(E,\Theta)}(\ispk)$.
Since there is an equivalence
$\widetilde{P}C^{\bullet}(\Theta',\Theta,UX)\simeq 
C^{\bullet}(R,\Theta,UX)$
of cosimplicial objects,
we see that $\widetilde{P}X\simeq PX$
for any $X$.
This completes the proof. 
\qqq

Let $\mathcal{T}$ be the full subcategory of $\mod_A(\ispk)$
consisting of $X$ such that
the unit map $X\to P{\rm Coex}(X)$ is an equivalence
\[ \mathcal{T}=\{X\in\mod_A(\ispk)|\ X\stackrel{\simeq}{\to}P{\rm Coex}(X)\}.\]
In the same way as in 
Proposition~\ref{prop:fundamental-proposition-module-embedding},
we obtain the following proposition.

\begin{proposition}
\label{prop:quasi-category-module-embedding}
The restriction of ${\rm Coex}$ to $\mathcal{T}$
is a fully faithful functor of
quasi-categories.
\end{proposition}

\subsection{Examples}

In this subsection we give some examples
of embeddings into quasi-categories of comodules. 

\begin{enumerate}

\item 
{({cf.~\cite[6.1.2]{Hess}})}
Let $k=S$ and let
$\psi: S\to MU$ be the unit map,
where $MU$ is the complex cobordism spectrum.
We have an adjunction
$MU\wedge(-):\isp\rightleftarrows
   \mod_{MU}(\isp):\psi^*$.
We denote by $MU\wedge MU$
the comonad on $\mod_{MU}(\isp)$ associated to
the adjoint pair $(MU\wedge(-),\psi^*)$.
If $X\in\isp$ is connective,  
then the map
$X\to \lim MU^{\bullet+1}X$ is an equivalence
by \cite[Thm.~6.5]{Bousfield}.
Hence $X\in\mathcal{T}$.
Let ${\rm Sp}^{\rm \ge 0}$ be the full subcategory
of $\isp$ consisting of connective spectra.
By Proposition~\ref{prop:quasi-category-module-embedding},
the functor
\[ MU\wedge(-): {\rm Sp}^{\ge 0}\longrightarrow
   \comod_{(MU,MU\wedge MU)}(\isp)\]
is fully faithful. 

\item 
Let $k=E(n)$ be the $n$th Johnson-Wilson spectrum 
at a prime $p$ and let
$\psi: S\to E(n)$ be the unit map.
We have an adjunction
$E(n)\wedge(-):\isp_{E(n)}\rightleftarrows
   \mod_{E(n)}(\isp_{E(n)}):\psi^*$.
We denote by $E(n)\wedge E(n)$
the comonad on $\mod_{E(n)}(\isp_{E(n)})$ associated to
the adjoint pair $(E(n)\wedge(-),\psi^*)$.
For any $X\in\isp_{E(n)}$,
the map
$X\to \lim E(n)^{\bullet+1}X$ is an equivalence
in $\isp_{E(n)}$
since any spectrum is $E(n)$-prenilpotent
by \cite[Thm.~5.3]{Hovey-Sadofsky}.
Hence the functor
\[ E(n)\wedge(-): \isp_{E(n)}\longrightarrow
   \comod_{(E(n),E(n)\wedge E(n))}(\isp_{E(n)})\]
is fully faithful.

\item
Let $k=K(n)$ and let
$\psi: S\to E_n$ be the unit map.
We have an adjunction
$L_{K(n)}(E_n\wedge(-)):\isp_{K(n)}\rightleftarrows
   \mod_{E_n}(\isp_{K(n)}):\psi^*$.
We denote by $L_{K(n)}(E_n\wedge E_n)$
the comonad on $\mod_{E_n}(\isp_{K(n)})$ associated to
the adjoint pair $(L_{K(n)}(E_n\wedge(-)),\psi^*)$.
For any $X\in\isp_{K(n)}$,
the map
$X\to \lim L_{K(n)}E_n^{\bullet+1}X$ is an equivalence
in $\isp_{K(n)}$
since any $K(n)$-local spectrum is $K(n)$-local $E_n$-nilpotent
by \cite[Prop.~A.3]{Devinatz-Hopkins}.
Hence the functor
\[ L_{K(n)}(E_n\wedge(-)): \isp_{K(n)}\longrightarrow
   \comod_{(E_n,L_{K(n)}(E_n\wedge E_n))}(\isp_{K(n)})\]
is fully faithful.

\item
Let $k=M(p)$ be the mod $p$ Moore spectrum
at a prime $p$. 
Let $\psi: S\to H\mathbb{F}_p$ be the unit map,
where $H\mathbb{F}_p$ is
the mod $p$ Eilenberg-Mac~Lane spectrum.
We have an adjunction
$H\mathbb{F}_p\wedge(-):\isp_{M(p)}\rightleftarrows
   \mod_{H\mathbb{F}_p}(\isp_{M(p)}):\psi^*$.
We denote by $H\mathbb{F}_p\wedge H\mathbb{F}_p$
the comonad on $\mod_{H\mathbb{F}_p}(\isp_{M(p)})$ associated to
the adjoint pair $(H\mathbb{F}_p\wedge(-),\psi^*)$.
If $X\in\isp_{M(p)}$ is connective,
then the map
$X\to \lim M\mathbb{F}_p^{\bullet+1}X$ is an equivalence
in $\isp_{M(p)}$
by \cite[Thm.~6.5]{Bousfield}, and 
hence $X\in\mathcal{T}$.
The functor
\[ H\mathbb{F}_p\wedge(-): {\rm Sp}_{M(p)}^{\ge 0}\longrightarrow
   \comod_{(H\mathbb{F}_p,H\mathbb{F}_p\wedge H\mathbb{F}_p)}(\isp_{M(p)})\]
is fully faithful,
where ${\rm Sp}_{M(p)}^{\ge 0}$ is the full subcategory
of ${\rm Sp}_{M(p)}$ consisting of connective spectra.


\end{enumerate}

\section{Quasi-category of discrete $G$-spectra}
\label{section:quasi-category-spG}

In this section we discuss
the underlying quasi-categories of 
modules and algebras over $\sp(G)_k$.
We show that the formulations of embeddings
of module categories in 
\S\ref{section:embedding_model_module} and 
\S\ref{section:embedding-quasi-categories}
are equivalent under some conditions.

\subsection{The quasi-category $\isp(G)_k$
as a comodule category}
\label{subsection: ispG-as-comodule}

By Proposition~\ref{proposition:k-local-U-V-adjunction},
we have a symmetric monoidal $\sp$-Quillen adjunction
$U: \sp(G)_k\rightleftarrows \spk:V$,
where $U$ is the forgetful functor 
and $V(-)={\rm Map}_c(G,-)$.
We denote by $\isp(G)_k$
the underlying quasi-category of $\sp(G)_k$.
The adjoint pair $(U,V)$ induces
an adjunction of quasi-categories
\[ U_k: \isp(G)_k\rightleftarrows\ispk:V_k.\]
We denote by $\Gamma$
the comonad on $\ispk$ associated to the
adjoint pair $(U_k,V_k)$.
We have the quasi-category of comodules 
\[ \comod_{\Gamma}(\ispk)\]
over $\Gamma$.
In this subsection we shall show that $U_k$ exhibits
$\isp(G)_k$ as comonadic over $\ispk$,
that is,
$\isp(G)_k$ is equivalent to 
$\comod_{\Gamma}(\ispk)$ 
under some conditions.

First, we consider the unlocalized version of the adjunction
$U:\isp(G)\rightleftarrows \isp:V$
and show that the forgetful functor
$U:\isp(G)\to \isp$
exhibits $\isp(G)$ as comonadic over $\isp$.
For this purpose,
we shall apply
the quasi-categorical Barr-Beck theorem by
Lurie \cite[4.7.4]{Lurie2}.
In particular,
we have to show that $U$ preserves
the limit of any cosimplicial
object of $\sp(G)$ that is split in $\sp$.
Since the limit of a diagram in the underlying quasi-category
of a simplicial model category is represented
by the homotopy limit of the simplicial model 
category,
we recall the homotopy limit of
a diagram in a simplicial model category.

We use a model of homotopy limits 
in \cite[Ch.~18]{Hirschhorn}.
For a small category $\mathcal{C}$ and 
an object $\alpha\in\mathcal{C}$,
we denote by $(\mathcal{C}\downarrow\alpha)$
the category of objects of $\mathcal{C}$ over $\alpha$.
An object of $(\mathcal{C}\downarrow\alpha)$
is a pair $(\beta,\sigma)$
where $\beta$ is an object of $\mathcal{C}$
and $\sigma$ is a map $\beta\to \alpha$ in $\mathcal{C}$.
A morphism from $(\beta,\sigma)$ to $(\beta',\sigma')$
in $(\mathcal{C}\downarrow\alpha)$ is a map
$\tau:\beta\to\beta'$ in $\mathcal{C}$
such that $\sigma=\sigma'\tau$. 
A map $\sigma:\alpha\to\alpha'$ in $\mathcal{C}$
induces a functor
$\sigma_*: (\mathcal{C}\downarrow\alpha)\to
(\mathcal{C}\downarrow\alpha')$
by composing with $\sigma$.
For a small category $\mathcal{C}$,
we denote by $B\mathcal{C}$ its classifying space (simplicial set),
which is given by applying the nerve functor to $\mathcal{C}$.
Let $\mathbf{M}$ be a simplicial model category.
For an object $\mathbf{Y}\in\mathbf{M}$
and a simplicial set $K$,
we denote by $\mathbf{Y}^K$
the power of $\mathbf{Y}$ by $K$.
For a functor $\mathbf{X}: \mathcal{C}\to \mathbf{M}$
such that $\mathbf{X}_{\alpha}$ is fibrant
for all $\alpha\in\mathcal{C}$,
the homotopy limit ${\rm holim}_{\mathcal{C}}^{\mathbf{M}}\,\mathbf{X}$ 
is defined to be
the equalizer of the maps
\[ \prod_{\alpha\in\mathcal{C}}
   (\mathbf{X}_{\alpha})^{B(\mathcal{C\downarrow\alpha})}\quad
   \stackrel{\phi}{\subrel{\psi}{\rightrightarrows}}\quad
   \prod_{(\sigma:\alpha\to\alpha')\in\mathcal{C}}
   (\mathbf{X}_{\alpha'})^{B(\mathcal{C}\downarrow\alpha)},\]
where the projection of the map $\phi$ 
on the factor $\sigma:\alpha\to\alpha'$
is the composition of the projection on the factor $\alpha$
with the map 
$(\mathbf{X}_{\alpha})^{B(\mathcal{C}\downarrow\alpha)}\to
 (\mathbf{X}_{\alpha'})^{B(\mathcal{C}\downarrow\alpha)}$
induced by the map
$\mathbf{X}_{\sigma}:\mathbf{X}_{\alpha}\to\mathbf{X}_{\alpha'}$
and the projection of the map $\psi$ on the 
factor $\sigma:\alpha\to\alpha'$ is the composition
of the projection on the factor $\alpha'$
with the map 
$(\mathbf{X}_{\alpha'})^{B(\mathcal{C}\downarrow\alpha')}\to
 (\mathbf{X}_{\alpha'})^{B(\mathcal{C}\downarrow\alpha)}$
induced by
the map $B(\sigma_*): B(\mathcal{C}\downarrow\alpha)\to
B(\mathcal{C}\downarrow\alpha')$.

Next,
we would like to
compare the homotopy limit of a diagram
in $\sp(G)$ with that of the diagram in $\sp$
obtained by applying the forgetful functor $U$. 
For this purpose,
we shall describe the limits and powers in $\sp(G)$  
in terms of those in $\sp$.

For a small category $\mathcal{C}$,
we let $\mathbf{X}: \mathcal{C}\to \sp(G)$
be a functor. 
We shall describe the limit
${\rm lim}_{\mathcal{C}}^{\sp(G)}\,\mathbf{X}$
in $\sp(G)$.
We have the induced functor $U\mathbf{X}:\mathcal{C}\to \sp$
and the limit
${\rm lim}_{\mathcal{C}}^{\sp}\, U\mathbf{X}$
in $\sp$.
By the functoriality of limit,
we can regard ${\rm lim}_{\mathcal{C}}^{\sp}\, U\mathbf{X}$
as an object in $\sp(G^{\delta})$.
We also have the induced functor
$\mathbf{X}^{\delta}:\mathcal{C}\to \sp(G^{\delta})$.
We see that ${\rm lim}_{\mathcal{C}}^{\sp}\, U\mathbf{X}$
can be identified with the limit of $\mathbf{X}^{\delta}$
in $\sp(G^{\delta})$:
${\rm lim}_{\mathcal{C}}^{\sp(G^{\delta})}\,\mathbf{X}^{\delta}\cong
{\rm lim}_{\mathcal{C}}^{\sp}\,U\mathbf{X}$.
Since $d:\sp(G^{\delta})\to \sp(G)$
is right adjoint,
the functor $d$ preserves limits,
and hence we obtain that
the limit ${\rm lim}_{\mathcal{C}}^{\sp(G)}\,\mathbf{X}$
in $\sp(G)$ is isomorphic to
$d({\rm lim}_{\mathcal{C}}^{\sp(G^{\delta})}\,\mathbf{X}^{\delta})$
and hence we obtain an isomorphism
\[ {\rm lim}_{\mathcal{C}}^{\sp(G)}\,\mathbf{X}\cong
   d({\rm lim}_{\mathcal{C}}^{\sp}\, U\mathbf{X}).\]

Let $X\in\sp(G)$ be a discrete symmetric $G$-spectrum
and $K$ a simplicial set. 
We shall describe the power $X^K$ of $X$ by $K$ in $\sp(G)$.
Note that the copower $X\otimes K$ 
is given by $X\wedge\Sigma^{\infty}_+ K$,
where $\Sigma^{\infty}_+K$
is the suspension spectrum of $K$
with the disjoint base point.
This implies that the functor $\delta$
preserves copower, that is,
there is a natural isomorphism
$(X\otimes K)^{\delta}\cong X^{\delta}\otimes K$
in $\sp(G^{\delta})$
for any $X\in\sp(G)$
and any simplicial set $K$.
By using the adjunction
of power and copower,
we see that
the functor $d$ preserves powers, 
that is,
$d(Y^K)\cong (dY)^K$
for any $Y\in \sp(G^{\delta})$
and any simplicial set $K$.
By the functoriality of power,
we can regard 
the power $(UX)^K$ of the symmetric spectrum $UX\in\sp$
by $K$ as an object of $\sp(G^{\delta})$
and identify
the power $(X^{\delta})^K$ in $\sp(G^{\delta})$
with $(UX)^K$.
Using $X\cong d(X^{\delta})$,
we obtain an isomorphism
\[ X^K\cong d((UX)^K)\]
in $\sp(G)$
for any $X\in \sp(G)$ and any simplicial set $K$.

Using these descriptions of limits
and powers in $\sp(G)$,
we shall describe the homotopy limit
of a diagram in $\sp(G)$. 
Let $\mathcal{C}$ be a small category
and let $\mathbf{X}:\mathcal{C}\to\sp(G)$
be a functor such that $\mathbf{X}_{\alpha}$
is fibrant for all $\alpha\in\mathcal{C}$.
By the functoriality of the homotopy limit,
we can regard
the homotopy limit
${\rm holim}_{\mathcal{C}}^{\sp}\, U\mathbf{X}$
in $\sp$ of the induced functor $U\mathbf{X}$
as an object of $\sp(G^{\delta})$.
Note that $U\mathbf{X}_{\alpha}$ is fibrant
for all $\alpha$ since $U$ preserves
fibrant objects by \cite[Cor.~5.3.3]{Behrens-Davis}.
Since $d$ preserves limits and powers,
we obtain an isomorphism
\[ {\rm holim}_{\mathcal{C}}^{\sp(G)}\, \mathbf{X}\cong
   d({\rm holim}_{\mathcal{C}}^{\sp}\, U\mathbf{X}).\] 

We use the following notation
for simplicity.
For a cosimplicial object $Y^{\bullet}$ in $\sp$
such that $Y^r$ is fibrant for all $r\ge 0$,
we denote by
\[ {\rm holim}_{\Delta}\, Y^{\bullet} \]
the homotopy limit 
${\rm holim}_{\Delta}^{\sp}\,Y^{\bullet}$ in $\sp$,
and for a cosimplicial object $Z^{\bullet}$ in $\sp(G)$
such that $Z^r$ is fibrant for all $r\ge 0$,
we denote by 
\[ {\rm holim}^G_{\Delta}\, Z^{\bullet}\]  
the homotopy limit of
$Z^{\bullet}$ in $\sp(G)$.
By the above argument,
we have an isomorphism
\[ {\rm holim}^G_{\Delta}\, Z^{\bullet}\cong
   d({\rm holim}_{\Delta}\,UZ^{\bullet}) \]
for any cosimplicial object $Z^{\bullet}$
in $\sp(G)$ such that
$Z^r$ is fibrant for all $r\ge 0$.

We have a canonical map
\[ U\,{\rm holim}^G_{\Delta}\ Z^{\bullet}\longrightarrow
   {\rm holim}_{\Delta}\, UZ^{\bullet}\]
in $\sp$.
We give a sufficient condition
to ensure that the canonical map
is an equivalence. 
The following lemma will be used 
to show that the forgetful functor $U$
preserves th limit of 
any $U$-split cosimplicial object of $\isp(G)$.

\begin{lemma}
\label{lemma:U-preserves-limit}
Let $Z^{\bullet}$ be a cosimplicial object in $\sp(G)$
such that $Z^r$ is fibrant for all $r\ge 0$.
We set $Z={\rm holim}^G_{\Delta}\ Z^{\bullet}$
and $W={\rm holim}_{\Delta}\, UZ^{\bullet}$.
If $G$ has finite virtual cohomological dimension,
the coaugmentation map $W\to UZ^{\bullet}$
induces an isomorphism 
$\pi_*(W)\stackrel{\cong}{\to}\pi^0\pi_*(UZ^{\bullet})$,
and $\pi^s\pi_*(UZ^{\bullet})=0$ for all $s>0$,
then the canonical map 
$UZ{\to} W$ is a weak equivalence.
\end{lemma}

\proof
Since $G$ has finite virtual cohomological dimension,
as in the second paragraph of the proof of 
\cite[Thm.~7.4]{Davis},
we can choose a positive integer $m$ and 
a fundamental neighborhood system $\mathcal{N}$ of 
the identity element of $G$
consisting of open normal subgroups $N$ 
such that the cohomological dimension ${\rm cd}(N)\le m$.
In the following of this proof,
we fix such a system $\mathcal{N}$.
Since $\mathcal{N}$ is a fundamental neighborhood system
of the identity,
we have 
$d(X)\cong {\rm colim}_{N\in\mathcal{N}}\,X^N$
for any $X\in\sp(G)$.
Hence we obtain an isomorphism
\[ \begin{array}{rcl}
   {\rm holim}^G_{\Delta}\, Z^{\bullet}&\cong&
   d({\rm holim}_{\Delta}\, UZ^{\bullet})\\[2mm]
   &\cong& {\rm colim}_{N\in\mathcal{N}}\,{\rm
       holim}_{\Delta}\, (Z^{\bullet})^N.
   \end{array}\]


Since the discrete $G$-spectrum
$Z^r$ is fibrant in $\sp(G)$,
by \cite[Prop.~3.3.1(2)]{Behrens-Davis},
$Z^r$ is fibrant in $\sp(N)$
for any $r\ge 0$ and any $N\in\mathcal{N}$.
Hence we see that   
the fixed points spectrum
$(Z^r)^N$ is equivalent to the homotopy fixed
points spectrum $(Z^r)^{hN}$ for any $r\ge 0$
and any $N\in \mathcal{N}$.
Since $N\in\mathcal{N}$ has finite cohomological dimension,
by \cite[Thm.~7.4]{Davis},
we have
$(Z^r)^N\simeq\, {\rm holim}_{\Delta}\, {\rm Map}_c(G^{\bullet+1},Z^r)^N$.
Hence we obtain
\[ {\rm holim}^G_{\Delta}\ Z^{\bullet}
    \simeq\,
   {\rm colim}_{N\in\mathcal{N}}\, {\rm holim}_{\Delta}\,{\rm holim}_{\Delta}\,
   {\rm Map}_c(G^{\bullet+1},Z^{\bullet})^N. \]

We fix $N\in\mathcal{N}$ and $k\ge 0$,
and consider the cosimplicial spectrum
${\rm Map}_c(G^{k+1},Z^{\bullet})^N$.
The Bousfield-Kan spectral sequence 
abutting to the homotopy groups
of ${\rm holim}_{\Delta}{\rm
  Map}_c(G^{k+1},Z^{\bullet})^N$
has the form
\[ {}_IE_2^{s,t}\cong
   \pi^s\pi_t{\rm Map}_c(G^{k+1},Z^{\bullet})^N
   \Longrightarrow
   \pi_{t-s}{\rm holim}_{\Delta}{\rm Map}_c(G^{k+1},Z^{\bullet})^N. \]
Since $\pi_*{\rm Map}_c(G^{k+1},Z^r)^N\cong
{\rm Map}_c(G^{k+1},\pi_*(Z^r))^N$ for any $r\ge 0$,
we see that 
${}_IE_2^{s,t}=0$ for $s>0$
and ${}_IE_2^{0,t}\cong {\rm Map}_c(G^{k+1},\pi_t(W))^N$
by the assumptions
on the cohomotopy groups
of $\pi_*(UZ^{\bullet})$
and the fact that 
${\rm Map}_c(G^{k+1},-)^N$
is an exact functor
from the category of discrete $G$-modules
to the category of abelian groups.
Note that $\pi_t(W)$ is a discrete $G$-module
since it is a submodule of the discrete $G$-module
$\pi_t(Z^0)$. 
Hence the spectral sequence collapses from the $E_2$-page
and we obtain an isomorphism
\[ \pi_*{\rm holim}_{\Delta}{\rm Map}_c(G^{k+1},Z^{\bullet})^N
   \cong {\rm Map}_c(G^{k+1},\pi_*(W))^N.\]

Now we consider
the cosimplicial spectrum
${\rm holim}_{[r]\in \Delta}
   {\rm Map}_c(G^{\bullet+1},Z^r)^N$.
The Bousfield-Kan spectral sequence
abutting to the homotopy groups
of ${\rm holim}_{\Delta}{\rm holim}_{[r]\in \Delta}
   {\rm Map}_c(G^{\bullet+1},Z^r)^N$
has the form
\begin{align}\label{align:spectral-sequence-type-II}
 {}_{II}E_2^{s,t}\cong
   \pi^s\pi_t{\rm holim}_{[r]\in \Delta}
   {\rm Map}_c(G^{\bullet+1},Z^r)^N
   \Longrightarrow
   \pi_{t-s}{\rm holim}_{\Delta}\,{\rm holim}_{[r]\in \Delta}\,
   {\rm Map}_c(G^{\bullet+1},Z^r)^N.
\end{align}
Since $\pi_t{\rm holim}_{[r]\in\Delta}{\rm Map}_c(G^{\bullet+1},Z^r)^N
\cong {\rm Map}_c(G^{\bullet+1},\pi_t(W))^N$,
we see that
${}_{II}E_2^{s,t}\cong H_c^s(N;\pi_t(W))$.

Since the cohomological dimension
${\rm cd}(N)$ is uniformly bounded for $N\in \mathcal{N}$,
by taking the colimit over $N\in\mathcal{N}$
of the spectral sequences 
(\ref{align:spectral-sequence-type-II}),
we obtain a spectral sequence
\[ {}_{III}E_2^{p,q}\cong\,{\rm colim}_{N\in\mathcal{N}}\,
   H_c^p(N;\pi_q(W))\Longrightarrow \pi_{q-p}(Z)\]
as in the proof of \cite[Thm.~7.4]{Davis}.
Since ${}_{III}E_2^{0,q}\cong\pi_q(W)$ and 
${}_{III}E_2^{p,q}=0$ for $p>0$,
we see that the inclusion map
$UZ\to W$ is a weak equivalence in $\sp$.
\qqq

Let $F:\mathcal{C}\to\mathcal{D}$ be
a functor of quasi-categories and
$X^{\bullet}$ a cosimplicial object of $\mathcal{C}$.
We say that $X^{\bullet}$ is $F$-split
if $FX^{\bullet}$ is a split cosimplicial object
of $\mathcal{D}$
(see \cite[4.7.3]{Lurie2} for
the definition of a split simplicial object).

We consider the forgetful functor $U:\isp(G)\to\isp$
and a $U$-split cosimplicial object of $\isp(G)$.
Note that $\isp(G)$
admits all small limits by \cite[Cor.~4.2.4.8]{Lurie}
since $\isp(G)$ is the underlying quasi-category
of the combinatorial simplicial model category
$\sp(G)$. 
In particular,
the limit of any cosimplicial object exists.
We obtain the following lemma
by Lemma~\ref{lemma:U-preserves-limit}.

\begin{lemma}
\label{lemma:preserve-split-U-cosimplicial-quasi-category}
If $G$ has finite virtual cohomological dimension,
then the forgetful functor $U:\isp(G)\to\isp$
preserves the limit of any $U$-split cosimplicial object
of $\isp(G)$.
\end{lemma}

\proof
We recall that $\mathbf{M}^{\circ}$
is the full simplicial subcategory of a simplicial model category 
$\mathbf{M}$ consisting of objects
that are both fibrant and cofibrant, 
and that $N(\mathbf{M}^{\circ})$
is the underlying quasi-category of $\mathbf{M}$,
where $N(-)$ is the simplicial nerve functor.

Let $X^{\bullet}$ be a cosimplicial object in $\isp(G)$
that is $U$-split.
By \cite[Prop.~4.2.4.4]{Lurie},
there is a cosimplicial object
$Y^{\bullet}$ in $\sp(G)^{\circ}$ such that
$N(Y^{\bullet})\simeq X^{\bullet}$.
Note that $UY^{\bullet}$ is a cosimplicial
object in $\sp^{\circ}$
since $U$ is a left Quillen functor
and preserves fibrant objects
by \cite[Cor.~5.3.3]{Behrens-Davis}.
In order to prove the lemma,
we have to show that the map
\[ \varphi: U({\rm holim}^G_{\Delta}Y^{\bullet})\longrightarrow
         {\rm holim}_{\Delta}\,UY^{\bullet}\]
is a weak equivalence. 
 
Since $UX^{\bullet}$ is a split cosimplicial object
in $\isp$,
$\pi_*(UX^{\bullet})=\pi_*(UY^{\bullet})$
is a split cosimplicial object in the category of graded modules.
In particular,
the coaugmentation map
induces an isomorphism
$\pi_*({\rm holim}_{\Delta}UY^{\bullet})
\stackrel{\cong}{\to}\pi^0\pi_*(UY^{\bullet})$,
and
$\pi^s\pi_*(UY^{\bullet})=0$
for all $s>0$.
By Lemma~\ref{lemma:U-preserves-limit},
we see that 
the map $\varphi$
is a weak equivalence.
This completes the proof.
\qqq

Using Lemma~\ref{lemma:preserve-split-U-cosimplicial-quasi-category},
we obtain the following proposition
which shows that $\isp(G)$ is comonadic over $\isp$.

\begin{proposition}
\label{prop:ispgk-as-comodules-no-localization}
If $G$ has finite virtual cohomological dimension,
then the forgetful functor $U:\isp(G)\to\isp$ 
exhibits $\isp(G)$ as comonadic over $\isp$,
that is,
we have an equivalence of quasi-categories
\[ \isp(G)\stackrel{\simeq}{\longrightarrow}
   \comod_{\Gamma}(\isp).\]
\end{proposition}

\proof
We shall use the quasi-categorical Barr-Beck theorem by
Lurie \cite[4.7.4]{Lurie2}.
We have to show that $U$ is conservative,
$\isp(G)$ admits a limit 
for any $U$-split cosimplicial object, and 
the limit of any $U$-split cosimplicial object
is preserved by $U$. 

By the definition of the weak equivalences in $\sp(G)$, 
we see that the forgetful functor $U:\isp(G)\to\isp$
is conservative.
Since $\isp(G)$ is the underlying quasi-category
of the simplicial model category $\sp(G)$,
the quasi-category $\isp(G)$ admits all small limits  
by \cite[Cor.~4.2.4.8]{Lurie}.
By Lemma~\ref{lemma:preserve-split-U-cosimplicial-quasi-category},
the limit of any $U$-split cosimplicial object of $\isp(G)$
is preserved by $U$.
This completes the proof.
\qqq

The following lemma is useful
to show that
other adjunctions are comonadic.

\begin{lemma}
\label{lemma:preservation-comonadicity}
Suppose we have a commutative diagram of quasi-categories
\[ \begin{array}{ccc}
     \mathcal{C}' & \stackrel{F'}{\longrightarrow} & 
     \mathcal{D}'\\
     \mbox{$\scriptstyle p$}\bigg\downarrow
     \phantom{\mbox{$\scriptstyle p$}} &&
     \phantom{\mbox{$\scriptstyle q$}}\bigg\downarrow
     \mbox{$\scriptstyle q$}\\[2mm]
     \mathcal{C} & \stackrel{F}{\longrightarrow}& 
     \mathcal{D}.
   \end{array}\]
We assume that $F$ and $F'$ are left adjoint functors,
and that $p$ and $q$ are conservative.
Furthermore,
we assume that, for any cosimplicial object $X^{\bullet}$
in $\mathcal{C}'$, if $p(X^{\bullet})$ admits a limit in 
$\mathcal{C}$,
then $X^{\bullet}$ admits a limit in $\mathcal{C}'$ and
that limit is preserved by $p$.
If $F$ exhibits $\mathcal{C}$ as comonadic over $\mathcal{D}$,
then $F'$ exhibits $\mathcal{C}'$ as comonadic over $\mathcal{D}'$. 
\end{lemma}

\proof
We shall use the quasi-categorical Barr-Beck theorem by
Lurie \cite[4.7.4]{Lurie2}.
We have to show that $F'$ is conservative,
$\mathcal{C}'$ admits a limit 
for any $F'$-split cosimplicial object, and 
the limit of any $F'$-split cosimplicial
object of $\mathcal{C}'$ is preserved by $F'$. 

First, we show that $F'$ is conservative.
Since $F$ exhibits $\mathcal{C}$ as comonadic over $\mathcal{D}$,
$F$ is conservative.
Combining this with the conservativeness of $p$,
we see that $F'$ is conservative as well.

Next, we let $X^{\bullet}$ be an $F'$-split cosimplicial
object of $\mathcal{C}'$.
By applying $p$,
we obtain an $F$-split cosimplicial object $p(X^{\bullet})$
in $\mathcal{C}$ since split cosimplicial objects
are preserved by any functor.
Since $F$ exhibits $\mathcal{C}$ as comonadic over $\mathcal{D}$,
$p(X^{\bullet})$ admits a limit and that limit is preserved by $F$.
By the assumption,
$X^{\bullet}$ admits a limit and the limit is preserved by $p$.
We see that the limit of $X^{\bullet}$ is preserved by $F'$
since the composition $Fp$ preserves the limit
and $q$ is conservative.
\qqq

Next, we consider the localized version
of the adjunction
$U_k:\isp(G)_k\rightleftarrows \ispk:V_k$
and would like to show that
the forgetful functor $U_k:\isp(G)_k\to \ispk$
exhibits $\isp(G)_k$ as comonadic over $\ispk$.
For this purpose,
we consider the following assumption
on the localization functor $L_k$.

\begin{assumption}[cf.~{\cite[Assumption~1.0.3]{Behrens-Davis}}]\rm
\label{condition:Behrens-Davis}
The localization functor $L_k$ on
the stable homotopy category of 
(non-equivariant) spectra
is given as a composite of
two localization functors $L_ML_T$,
where $L_T$ is a smashing localization and 
$L_M$ is a localization with respect to a finite spectrum $M$. 
\end{assumption}

We note that the $K(n)$-localization
$L_{K(n)}$ satisfies Assumption~\ref{condition:Behrens-Davis},
where $K(n)$ is the $n$th Morava $K$-theory 
at a prime $p$.
Let $E(n)$ be the $n$th Johnson-Wilson spectrum 
and $F(n)$ a finite spectrum of type $n$ at $p$.
The $E(n)$-localization $L_{E(n)}$ is smashing
by \cite[Thm.~7.5.6]{Ravenel}, and
the $K(n)$-localization is given as 
the composite $L_{F(n)}L_{E(n)}$ 
(see, for example, \cite[Prop.~7.10]{Hovey-Strickland}). 

In order to compare the comonadicity
of $\sp(G)_k$ and $\sp(G)$,
we consider the functor
\[ M\wedge(-): \isp(G)_k\to\isp(G) \]
given by smashing with a (non-equivariant) finite spectrum $M$.
The functor $M\wedge (-)$ is right adjoint
to the functor $(DM\wedge(-))_k$,
where $DM$ is the $S$-dual of $M$
and $(-)_k$ is the $k$-localization functor
on $\isp(G)$.
We shall show that 
$M\wedge(-)$ is conservative
if the localization functor $L_k$
satisfies $L_ML_k\simeq L_k$.

\begin{lemma}
\label{lemma:conservativeness-M-smashing}
If the localization functor $L_k$ satisfies
$L_ML_k\simeq L_k$ for a finite spectrum $M$,
then the functor $M\wedge(-): \isp(G)_k\to \isp(G)$
is conservative.
\end{lemma}

\proof
Suppose $f: X\to Y$ is a map in $\isp(G)_k$
such that $M\wedge f$ is an equivalence in $\isp(G)$.
We have to show that $Uf$ is a $k$-equivalence,
where $U:\isp(G)_k\to\isp$ is
the forgetful functor.
Since $M$ is finite and $M\wedge Uf$ is an equivalence,
$M\wedge L_kUf\simeq L_k(M\wedge Uf)$ 
is also an equivalence.
Hence $L_ML_kUf\simeq L_kUf$ is an equivalence.
This completes the proof.
\qqq

The following theorem shows that
$\isp(G)_k$ is comonadic over $\isp_k$
under some conditions. 

\begin{theorem}
\label{thm:ispgk-as-comodules-localized-version}
If $G$ has finite virtual cohomological dimension
and the localization functor $L_k$ satisfies
Assumption~\ref{condition:Behrens-Davis},
then the forgetful functor $U:\isp(G)_k\to\ispk$
exhibits $\isp(G)_k$ as comonadic over $\ispk$,
that is,
we have an equivalence of quasi-categories
\[ \isp(G)_k\stackrel{\simeq}{\longrightarrow}
   \comod_{\Gamma}(\ispk).\]
\end{theorem}

\proof
Suppose $L_k\simeq L_ML_T$,
where $M$ is a finite spectrum and $L_T$ is smashing.
We shall apply Lemma~\ref{lemma:preservation-comonadicity}
for the following diagram
\[ \begin{array}{ccc}
    \isp(G)_k & \stackrel{U_k}
    {\hbox to 10mm{\rightarrowfill}} & \ispk\\[0.5mm]
    \mbox{$\scriptstyle M\wedge(-)$}\bigg\downarrow
    \phantom{\mbox{$\scriptstyle M\wedge(-)$}} & & 
    \phantom{\mbox{$\scriptstyle M\wedge(-)$}}
    \bigg\downarrow 
    \mbox{$\scriptstyle M\wedge(-)$}\\[2mm]
    \isp(G) & \stackrel{U}
    {\hbox to 10mm{\rightarrowfill}} & \isp,\\ 
   \end{array}\]
where the vertical arrows are given by
smashing with the finite spectrum $M$.

First, we have to show that
the diagram is commutative.
Let $X$ be an object of $\isp(G)_k$
which is represented by a fibrant and cofibrant object $Y$ 
in $\sp(G)_k$.
Since $U_kX$ in $\ispk$
is represented by $L_kUY$ in $\sp$,
we have to show that 
the natural map $M\wedge UY\to M\wedge L_kUY$
is a weak equivalence in $\sp$.
This follows from the assumption
that $L_k\simeq L_ML_T$
and the fact that $UY$ is $T$-local
by \cite[Prop.~6.1.7(2)]{Behrens-Davis}.
 
The horizontal arrows are left adjoint functors,
and the vertical arrows are conservative
by Lemma~\ref{lemma:conservativeness-M-smashing}.
Since $\isp(G)_k$ is the underlying quasi-category
of the combinatorial simplicial model category $\sp(G)_k$,
it admits all small limits
by \cite[Cor.~4.2.4.8]{Lurie}.
Since $M\wedge (-): \isp(G)_k\to\ispk$
is right adjoint,
it preserves all limits
by \cite[Prop.~5.2.3.5]{Lurie}.
By Proposition~\ref{prop:ispgk-as-comodules-no-localization},
$U$ exhibits $\isp(G)$ as comonadic over $\isp$.
Hence, by Lemma~\ref{lemma:preservation-comonadicity},
$U_k$ exhibits $\isp(G)_k$ as comonadic over $\ispk$.
\qqq

\subsection{The quasi-category of algebra objects in $\isp(G)_k$}

In this subsection we shall compare
the quasi-category of algebra objects
in $\isp(G)_k$ 
with the underlying quasi-category of the simplicial
model category of monoid objects in $\sp(G)_k$.
Furthermore, we shall show
that the quasi-category of algebra objects
in $\isp(G)_k$ is comonadic over 
the quasi-category of algebra objects in $\ispk$.

First, we recall a model structure
on the category of monoid objects in $\sp(G)_k$.
We denote by ${\rm Alg}(\sp(G)_k)$ the category 
of monoid objects in $\sp(G)_k$, and
let $F: {\rm Alg}(\sp(G)_k)\to \sp(G)_k$
be the forgetful functor.
By \cite[Thm.~4.1(3)]{Schwede-Shipley},
${\rm Alg}(\sp(G)_k)$ supports a model structure
as follows.
A map $f:X\to Y$ in ${\rm Alg}(\sp(G)_k)$
is said to be 
\begin{itemize}
\item
a weak equivalence if $F(f)$ is a 
weak equivalence in $\sp(G)_k$,
\item
a fibration if $F(f)$ is a 
fibration in $\sp(G)_k$, and
\item
a cofibration if  
it has the right lifting property
with respect to all maps which are both
fibrations and weak equivalences.
\end{itemize}
By Proposition~\ref{prop:spGk-combinatorial},
Theorem~\ref{thm:spgk-symmetric-monoidal-model-category}
and Proposition~\ref{prop:spgk-satisfies-monoid-axiom},
$\sp(G)_k$ is a combinatorial 
symmetric monoidal simplicial model category
which satisfies the monoid axiom.
We see that 
${\rm Alg}(\sp(G)_k)$ is a simplicial model category
and 
the forgetful functor $F:{\rm Alg}(\sp(G)_k)\to\sp(G)_k$
is a simplicial right Quillen functor
by \cite[4.1.4]{Lurie2}. 

We compare the quasi-category of algebra
objects in $\isp(G)_k$
with the underlying quasi-category of
the simplicial model category
${\rm Alg}(\sp(G)_k)$.
Let ${\rm Alg}(\isp(G)_k)$ be the quasi-category
of algebra objects in $\isp(G)_k$.
By \cite[1.3.4 and 4.1.4]{Lurie2},
there is an equivalence of quasi-categories
\[ N({\rm Alg}(\sp(G)_k)^{\circ})\simeq {\rm Alg}(\isp(G)_k),\]
where ${\rm Alg}(\sp(G)_k)^{\circ}$
is the full simplicial subcategory of ${\rm Alg}(\sp(G)_k)$
consisting of objects that are both fibrant and cofibrant,
and $N(-)$ is the simplicial nerve functor.

The forgetful functor $U: \sp(G)_k\to \spk$
induces a functor $U: {\rm Alg}(\sp(G)_k)\to{\rm Alg}(\spk)$.
We construct a right adjoint to 
the functor $U:{\rm Alg}(\sp(G)_k)\to{\rm Alg}(\spk)$.
For $Y\in {\rm Alg}(\spk)$,
we have an object ${\rm Map}_c(G,Y)$ in $\sp(G)_k$. 
We consider a map
\[ {\rm Map}_c(G,Y)\wedge {\rm Map}_c(G,Y)
   \longrightarrow {\rm Map}_c(G,Y) \]
in $\sp(G)_k$, which is the adjoint to the map
\[ U({\rm Map}_c(G,Y)\wedge {\rm Map}_c(G,Y))\cong
   U({\rm Map}_c(G,Y))\wedge U({\rm Map}_c(G,Y))
   \stackrel{{\rm ev}(e)\wedge {\rm ev}(e)}
   {\hbox to 17mm{\rightarrowfill}}  Y\wedge Y
   \stackrel{m}{\hbox to 10mm{\rightarrowfill}}Y,\]
where ${\rm ev}(e)$ is the evaluation map at
the identity element $e\in G$
and $m$ is the multiplication map on $Y$.
By this map,
we can regard ${\rm Map}_c(G,Y)$
as an object in ${\rm Alg}(\sp(G)_k)$.
Hence we obtain a functor
\[ V: {\rm Alg}(\spk)\longrightarrow {\rm Alg}(\sp(G)_k) \]
given by $V(Y)={\rm Map}_c(G,Y)$.
We see that $V$ is a right adjoint to the forgetful
functor $U$. 

\begin{proposition}
\label{prop:Alg-UV-adjunction}
The adjoint pair of functors
\[ U: {\rm Alg}(\sp(G)_k)\rightleftarrows 
      {\rm Alg}(\spk):V\]
is a simplicial Quillen adjunction.
\end{proposition}

\proof
Let $F:{\rm Alg}(\sp(G)_k)\to \sp(G)_k$
be the forgetful functor.
We consider the following commutative diagram
\[ \begin{array}{ccc}
    {\rm Alg}(\spk)&\stackrel{V}{\longrightarrow}&
    {\rm Alg}(\sp(G)_k)\\[1mm]
    \mbox{$\scriptstyle F$}\bigg\downarrow
    \phantom{\mbox{$\scriptstyle F$}}
    & & 
    \phantom{\mbox{$\scriptstyle F$}}
    \bigg\downarrow\mbox{$\scriptstyle F$}\\[2mm]
    \spk & \stackrel{V}{\longrightarrow} &
    \sp(G)_k.
   \end{array}\]
By Proposition~\ref{proposition:k-local-U-V-adjunction},
$V: \spk\to\sp(G)_k$ is a right Quillen functor.
This implies that $V: {\rm Alg}(\spk)\to {\rm Alg}(\sp(G)_k)$
is also a right Quillen functor.
\qqq

By Proposition~\ref{prop:Alg-UV-adjunction},
we have an adjunction of quasi-categories
\[ U_k: {\rm Alg}(\isp(G)_k)\rightleftarrows 
      {\rm Alg}(\ispk):V_k.\]
This induces a map of quasi-categories
\[ {\rm Alg}(\isp(G)_k)\longrightarrow
   \comod_{\Gamma}({\rm Alg}(\ispk)),\]
where $\Gamma$ is the comonad on ${\rm Alg}(\ispk)$
associated to the adjunction $(U_k,V_k)$.
The following theorem 
shows that ${\rm Alg}(\isp(G)_k)$
is comonadic over ${\rm Alg}(\ispk)$.

\begin{theorem}
\label{theorem:rectification_AlgG}
Let $G$ be a profinite group that has finite virtual
cohomological dimension.
We assume that the localization functor
$L_k$ satisfies 
Assumption~\ref{condition:Behrens-Davis}.
Then the forgetful functor 
$U_k: {\rm Alg}(\isp(G)_k)\to {\rm Alg}(\ispk)$
exhibits ${\rm Alg}(\isp(G)_k)$
as comonadic over ${\rm Alg}(\ispk)$,
that is, we have an equivalence 
of quasi-categories
\[ {\rm Alg}(\isp(G)_k)\stackrel{\simeq}{\longrightarrow}
   \comod_{\Gamma}({\rm Alg}(\ispk)).\]
\end{theorem}

\proof
We shall apply Lemma~\ref{lemma:preservation-comonadicity}
for the commutative diagram
\[ \begin{array}{ccc}
    {\rm Alg}(\isp(G)_k)&
    \stackrel{U_k}{\longrightarrow}&
    {\rm Alg}(\ispk)\\[1mm]
    \mbox{$\scriptstyle F$}\bigg\downarrow
    \phantom{\mbox{$\scriptstyle F$}} & & 
    \phantom{\mbox{$\scriptstyle F$}}
    \bigg\downarrow\mbox{$\scriptstyle F$}\\[2mm]
    \isp(G)_k & \stackrel{U_k}{\longrightarrow}&
    \ispk.\\
   \end{array}\]
The horizontal arrows are left adjoint functors
and the vertical arrows are conservative.
Since ${\rm Alg}(\sp(G)_k)$
is a combinatorial simplicial model category
and ${\rm Alg}(\isp(G)_k)$ is its underlying
quasi-category,
${\rm Alg}(\isp(G)_k)$ admits all small limits 
by \cite[Cor.~4.2.4.8]{Lurie}.
Since the forgetful functor 
$F: {\rm Alg}(\isp(G)_k)\to \isp(G)_k$
is right adjoint,
it preserves limits 
by \cite[Prop.~5.2.3.5]{Lurie}.
By Theorem~\ref{thm:ispgk-as-comodules-localized-version},
$U_k: \isp(G)_k\to\ispk$ exhibits $\isp(G)_k$
as comonadic over $\ispk$.
Hence the theorem follows from 
Lemma~\ref{lemma:preservation-comonadicity}.
\qqq

\subsection{The quasi-category of 
module objects in $\isp(G)_k$}

In this subsection we shall
show that the quasi-category
of module objects in $\isp(G)_k$ is comonadic over
the quasi-category of module
objects in $\ispk$.  

First,
we compare the quasi-category of module objects
in $\isp(G)_k$
and the underlying quasi-category 
of the simplicial model category of module objects in $\sp(G)_k$. 
Let $B$ be a monoid object in $\sp(G)_k$.
We assume that $B$ is cofibrant in $\sp(G)_k$.
We denote by $UB$ the underlying monoid object
in $\spk$.
Note that $UB$ is cofibrant
in $\spk$
by Proposition~\ref{proposition:k-local-U-V-adjunction}.
By \cite[1.3.4 and 4.3.3]{Lurie2},
the underlying quasi-categories of $\mod_{B}(\sp(G)_k)$
and $\mod_{UB}(\spk)$ are
equivalent to $\mod_{B}(\isp(G)_k)$ and
$\mod_{UB}(\ispk)$, respectively.

We have the forgetful functor
$U: \mod_B(\sp(G)_k)\to\mod_{UB}(\spk)$.
We shall construct a right adjoint to $U$.
For $M\in\mod_{UB}(\spk)$,
we regard ${\rm Map}_c(G,M)$
as an object in $\sp(G)_k$.
We consider a map
\[ B\wedge {\rm Map}_c(G,M)\longrightarrow 
   {\rm Map}_c(G,M) \]
in $\sp(G)_k$, 
which is an adjoint to the map  
\[ U(B\wedge {\rm Map}_c(G,M))\cong
   UB\wedge U{\rm Map}_c(G,M)\stackrel{{\rm id}\wedge {\rm ev}(e)}
{\hbox to 15mm{\rightarrowfill}}
   UB\wedge M
   \stackrel{a}{\hbox to 10mm{\rightarrowfill}}
   M,\]
where ${\rm ev}(e)$ is the evaluation map
at the identity element $e\in G$ and
$a$ is the action map on $M$.
This defines a $B$-module structure
on ${\rm Map}_c(G,M)$
and we see that ${\rm Map}_c(G,M)$
is an object in $\mod_B(\sp(G)_k)$.
Hence we obtain a functor
\[ V: \mod_{UB}(\spk)\longrightarrow \mod_B(\sp(G)_k)\]
given by $V(M)={\rm Map}_c(G,M)$.
We see that $V$ is a right adjoint to
the forgetful functor $U$, and hence
we have an adjunction
\[ U: \mod_B(\sp(G)_k)\rightleftarrows\mod_{UB}(\spk):V.\]  

\begin{lemma}
The adjoint pair $(U,V)$ is a $\sp$-Quillen adjunction.
\end{lemma}

\proof
This follows from 
Proposition~\ref{proposition:k-local-U-V-adjunction}.
\qqq

The $\sp$-Quillen adjunction
$(U,V)$ induces an adjunction of quasi-categories
\[ U_k: \mod_B(\isp(G)_k)\rightleftarrows
      \mod_{UB}(\ispk):V_k.\]
Let $\Gamma$ be the comonad on $\mod_{UB}(\ispk)$
associate to the adjoint pair $(U_k,V_k)$,
and let
\[ \comod_{(UB,\Gamma)}(\ispk)=\comod_{\Gamma}(\mod_{UB}(\ispk))\]
be the quasi-category of comodules over $\Gamma$.
The following theorem shows that
$\mod_B(\isp(G)_k)$ is comonadic over 
$\mod_{UB}(\ispk)$ under some conditions.

\begin{theorem}
\label{theorem:equivalence-BG-Gamma}
Let $G$ be a profinite group that has finite virtual cohomological 
dimension.
We assume that the localization functor
$L_k$ satisfies 
Assumption~\ref{condition:Behrens-Davis}.
Then the forgetful functor
$U_k: \mod_B(\isp(G)_k)\to \mod_{UB}(\ispk)$
exhibits $\mod_B(\isp(G)_k)$
as comonadic over $\mod_{UB}(\ispk)$,
that is, 
we have an equivalence of quasi-categories
\[ \mod_B(\isp(G)_k)\stackrel{\simeq}{\longrightarrow}
  \comod_{(UB,\Gamma)}(\ispk).\]
\end{theorem}

\proof
The theorem follows in the same way
as the proof of Theorem~\ref{theorem:rectification_AlgG}
by applying 
Lemma~\ref{lemma:preservation-comonadicity}
for the commutative diagram
\[ \begin{array}{ccc}
    \mod_B(\isp(G)_k)&
    \stackrel{U_k}{\longrightarrow}&
    \mod_{UB}(\ispk)\\[1mm]
    \phantom{\mbox{$\scriptstyle F$}}\bigg\downarrow
    \phantom{\mbox{$\scriptstyle F$}} & & 
    \phantom{\mbox{$\scriptstyle F$}}
    \bigg\downarrow\phantom{\mbox{$\scriptstyle F$}}\\[2mm]
    \isp(G)_k & \stackrel{U_k}{\longrightarrow}&
    \ispk,\\
   \end{array}\]
where the vertical arrows are forgetful functors.
\qqq

\subsection{Equivalence of the two formulations}
\label{subsection:equivalence-two-formulations}

In this subsection we shall show that
Proposition~\ref{prop:fundamental-proposition-module-embedding} and 
Proposition~\ref{prop:quasi-category-module-embedding}
are equivalent
in some circumstances.

First, we recall 
the relationship between 
functors of $\sp$-model categories and
of the underlying quasi-categories.
The homotopy category of the underlying quasi-category
of a $\sp$-model category $\mathbf{M}$ is equivalent to
the homotopy category of $\mathbf{M}$
as ${\rm Ho}(\sp)$-enriched categories.
Let 
\[ F:\mathbf{M}\rightleftarrows \mathbf{N}:G \]
be a $\sp$-Quillen adjunction
between combinatorial $\sp$-model categories.
The adjunction $(F,G)$ induces an adjunction
\[ \mathcal{F}: \mathcal{M}\rightleftarrows
   \mathcal{N}:\mathcal{G} \]
of quasi-categories,
where $\mathcal{M}$ and $\mathcal{N}$
are the underlying quasi-categories
of $\mathbf{M}$ and $\mathbf{N}$,
respectively.
The adjunction $(\mathcal{F},\mathcal{G})$
induces an adjunction 
\[ {\rm Ho}(\mathcal{F}): {\rm Ho}(\mathcal{M})
\rightleftarrows {\rm Ho}(\mathcal{N}):{\rm Ho}{G} \]
on the homotopy categories.
The adjunction
$({\rm Ho}(\mathcal{F}),{\rm Ho}(\mathcal{G}))$
is identified with
the derived adjunction
$(\mathbb{L}F,\mathbb{R}G)$
under the equivalences 
${\rm Ho}(\mathbf{M})\simeq {\rm Ho}(\mathcal{M})$ 
and ${\rm Ho}(\mathbf{N})\simeq {\rm Ho}(\mathcal{N})$.

Since $\mathbf{M}$ is a $\sp$-model category,
$\mathbf{M}$ is a stable model category
by \cite[Lem.~3.5.2]{Schwede-Shipley2},
and hence ${\rm Ho}(\mathbf{M})$ is a triangulated category. 
Suppose that $\mathbf{T}$ is a triangulated subcategory
of ${\rm Ho}(\mathbf{M})$.
We have the corresponding full subcategory
$\mathcal{T}$ of $\mathcal{M}$.
The total left derived functor
$\mathbb{L}(F): {\rm Ho}(\mathbf{M})\to
{\rm Ho}(\mathbf{N})$ restricted to $\mathbf{T}$
is fully faithful as an ${\rm Ho}(\sp)$-enriched functor 
if and only if
$\mathcal{F}: \mathcal{M}\to \mathcal{N}$
restricted to $\mathcal{T}$ 
is fully faithful.

Let $A$ be a monoid object in $\sp$.
We regard $A$ as a monoid object in $\sp(G)$
with trivial $G$ action.
Let $\varphi: A\to B$ be a map of monoid objects
in $\sp(G)$.
We assume that $A$ is cofibrant in $\sp$
and that $UB$ is fibrant and cofibrant in $\sp$.
We have the $\sp$-Quillen adjunction
\[
{\rm Ex}:\mod_A(\spk)\rightleftarrows\mod_B(\sp(G)_k):{\rm Re}, 
\]
which induces the adjunction 
$\mathbb{L}{\rm Ex}:{\rm Ho}(\mod_A(\spk))\rightleftarrows
{\rm Ho}(\mod_B(\sp(G)_k)):\mathbb{R}{\rm Re}$
of the homotopy categories.
The $\sp$-Quillen adjunction $({\rm Ex},{\rm Re})$
also induces the adjunction of quasi-categories
\[
{\mathcal Ex}:\mod_A(\ispk)\rightleftarrows
          \mod_B(\isp(G)_k):{\mathcal Re}. 
\]
We can regard this adjunction 
as a lifting of the adjunction 
$(\mathbb{L}{\rm Ex}, \mathbb{R}{\rm Re})$.

On the other hand,
we have the adjunction of quasi-categories
\[ U_k: \mod_A(\ispk)\rightleftarrows
        \mod_{UB}(\ispk):V_k,\]
which induces the adjunction of quasi-categories
\[ {\rm Coex}: \mod_A(\ispk)\rightleftarrows
    {\rm Comod}_{(UB,\Theta)}(\ispk):P.\]    
In this subsection,
under some conditions,
we shall show that
${\rm Comod}_{(UB,\Theta)}(\ispk)$
is equivalent to $\mod_B(\isp(G)_k)$
and the functor ${\mathcal Ex}$ is equivalent to
${\rm Coex}$
under this equivalence.
This implies that the right adjoint ${\mathcal Re}$
is equivalent to $P$, 
and the full subcategory $\mathcal{T}$
of $\mod_A(\ispk)$ corresponds to the thick subcategory
$\mathbf{T}$ of ${\rm Ho}(\mod_A(\spk))$,
where $\mathcal{T}$ is the full subcategory
of $\mod_A(\ispk)$ consisting of objects $X$
such that the unit map $X\to P{\rm Coex}(X)$
is an equivalence and 
the thick subcategory
$\mathbf{T}$ of ${\rm Ho}(\mod_A(\spk))$
consists of objects $Y$
such that the unit map 
$Y\to \mathbb{R}{\rm Re}\mathbb{L}{\rm Ex}(Y)$
is an isomorphism in the homotopy category. 
This means that 
the formulation of embeddings 
of module categories
in terms of model categories 
in Proposition~\ref{prop:fundamental-proposition-module-embedding} 
and that in terms of quasi-categories  
in Proposition~\ref{prop:quasi-category-module-embedding}
are equivalent.

We begin with constructing a map
which compares the two-sided bar construction
with the functor ${\rm Map}_c(G,-)$.
For $M\in \mod_{UB}(\sp)$,
we denote by
\[ B(UB,A,M)=|B_{\bullet}(UB,A,M)| \]
the geometric realization of the bar construction
$B_{\bullet}(UB,A,M)$. 
We define a map
\[ \Psi_M: B(UB,A,M)\longrightarrow
         U{\rm Map}_c(G,M)\] 
by applying $U$ to the map
$B(B,A,M)\to {\rm Map}_c(G,M)$
in $\sp(G)$
that is adjoint to 
the map $UB(B,A,M)\cong B(UB,A,M)\to M$
induced by the action of $UB$ on $M$.
In particular, we have a map   
\[ \Psi_{UB}: B(UB,A,UB)\longrightarrow
         U{\rm Map}_c(G,UB).\]


We set $\psi=U\varphi: A\to UB$.
Recall that there is a Quillen adjunction
\[ UB\wedge_A(-): \mod_A(\spk)\rightleftarrows
   \mod_{UB}(\spk): \psi^*, \]
which induces an adjunction
\[ UB\wedge^{\mathbb{L}}_A(-):
   {\rm Ho}(\mod_A(\spk))\rightleftarrows
   {\rm Ho}(\mod_{UB}(\spk)):\mathbb{R}\psi^*\]
between the homotopy categories,
where $UB\wedge^{\mathbb{L}}_A(-)$
is the total left derived functor
of $UB\wedge_A(-)$ and
$\mathbb{R}\psi^*$
is the total right derived functor
of $\psi^*$.
The following lemma shows that
the bar construction $B(UB,A,-)$
is a model of the composition
$UB\wedge^{\mathbb{L}}_A(\mathbb{R}\psi^*(-))$
of the functors. 

\begin{lemma}
\label{lemma:bar-construction-representative-UB-A-M}
If $M$ is a fibrant and cofibrant object
in $\mod_{UB}(\spk)$,
then 
$B(UB,A,M)$ represents
$UB\wedge^{\mathbb{L}}_A (\mathbb{R}\psi^*M)$
in ${\rm Ho}(\mod_{UB}(\spk))$.
\end{lemma}

\proof
Since $M$ is fibrant in $\mod_{UB}(\spk)$, 
$M$ represents $\mathbb{R}\psi^*M$
in ${\rm Ho}(\mod_A(\spk))$.
If $Q_AM\to M$ is a cofibrant replacement
in $\mod_A(\spk)$,
then $UB\wedge_A^{\mathbb{L}}(\mathbb{R}\psi^*M)$
is represented by $UB\wedge_A Q_AM$.
By \cite[Lem.~4.1.9]{Shipley2},
we have an equivalence
$B(UB,A,Q_AM)\stackrel{\simeq}{\to}
UB\wedge_A Q_AM$.     

We shall show that 
there is an equivalence
$B(UB,A,Q_AM)\stackrel{\simeq}{\to}
B(UB,A,M)$.
For any $r\ge 0$,
$UB\wedge A^{\wedge r}$ is cofibrant in $\sp$.
This implies an equivalence
$UB\wedge A^{\wedge r} \wedge Q_AM\stackrel{\simeq}{\to}
UB\wedge A^{\wedge r} \wedge M$
by \cite[Lem.~5.4.4]{HSS}.
Hence we obtain an equivalence
$|B_\bullet(UB,A,Q_AM)|\stackrel{\simeq}{\to}
|B_\bullet(UB,A,M)|$ 
by \cite[Cor.~4.1.6]{Shipley2}.
\qqq

The following lemma shows that
it suffices to show that 
$\Psi_{UB}$ is a $k$-local equivalence
in order to ensure that $\Psi_M$ is a $k$-local equivalence.

\begin{lemma}
\label{lemma:psi-ub-imply-psi-m}
Let $M$
be a cofibrant and fibrant object in $\mod_{UB}(\spk)$.
If $\Psi_{UB}$ is a $k$-local equivalence,
then $\Psi_M$ is also a $k$-local equivalence.
\end{lemma}

\proof
We have an isomorphism
between $B(UB,A,M)$ and $B(UB,A,UB)\wedge_{UB}M$,
and an equivalence 
between $U{\rm Map}_c(G,M)$ and $U{\rm Map}_c(G,UB)\wedge_{UB}M$.
Since $M$ is cofibrant in $\mod_{UB}(\sp)$,
the $k$-local equivalence $\Psi_{UB}$ 
induces a $k$-local equivalence
\[ B(UB,A,UB)\wedge_{UB}M\stackrel{\simeq_k}
   {\longrightarrow} U{\rm Map}_c(G,UB)\wedge_{UB}M\] 
by \cite[Lem.~5.4.4]{HSS}.
This completes the proof.
\qqq

There is an adjunction of quasi-categories
\[ UB\wedge_A(-): \mod_A(\ispk)\rightleftarrows
   \mod_{UB}(\ispk): \psi^*,\] 
and hence we obtain a comonad $\Theta$ on $\mod_{UB}(\ispk)$
and a quasi-category of comodules
\[ \comod_{(UB,\Theta)}(\ispk)=\comod_{\Theta}(\mod_{UB}(\ispk)).\]

To ease notation,
we set
\[ \begin{array}{rcl}
    \mathcal{C}&=& \mod_{UB}(\ispk)^{\rm op},\\[2mm]
    \mathcal{C}(G)&=& \mod_B(\isp(G)_k)^{\rm op},\\[2mm]
    \mathcal{D}&=&\mod_A(\ispk)^{\rm op}.
   \end{array}\]
We have an adjunction of quasi-categories
$V: \mathcal{C}\rightleftarrows \mathcal{C}(G): U$.
By \cite[4.7.4]{Lurie2},
there is an endomorphism monad of $U$,
and hence  
we have a monad 
$\Gamma\in {\rm Alg}({\rm End}(\mathcal{C}))$ 
and a left $\Gamma$-module 
$\overline{U}\in {\rm Mod}_{\Gamma}
({\rm Fun}(\mathcal{C}(G),\mathcal{C}))$.
Note that $\Gamma$ is a lifting of $UV$
and $\overline{U}$ is a lifting of $U$.

We set 
\[ \begin{array}{rcl}
     H&=&UB\wedge_A(-): \mathcal{D}\to\mathcal{C},\\[2mm]
     H'&=&B\wedge_A(-):\mathcal{D}\to \mathcal{C}(G),\\[2mm]
     F&=&\psi^*:\mathcal{C}\to\mathcal{D}.\\
   \end{array}\]
Note that $H=UH'$ is the right adjoint to $F$.
Hence there is an endomorphism monad of $H$,
which consists of a monad $\Theta\in {\rm Alg}({\rm End}(C))$
together with $\overline{H}\in 
{\rm Mod}_{\Theta}({\rm Fun}(\mathcal{D},\mathcal{C}))$
that is a lifting of $H$. 
We note that the functor $H=UH'$ lifts to a left $\Gamma$-module
$\overline{U}H'\in 
{\rm Mod}_{\Gamma}({\rm Fun}(\mathcal{D},\mathcal{C}))$.

We would like to show that
the monad $\Gamma$ together with
the left $\Gamma$-module object $\overline{U}H'$
is an endomorphism monad of $H$.
For this purpose, 
we consider the composite map
\[ \Gamma\stackrel{{\rm id}_{\Gamma}\times u}{\hbox to 10mm{\rightarrowfill}}
   \Gamma H F
   \stackrel{a\times {\rm id}_F}{\hbox to 10mm{\rightarrowfill}}
   H F,\]
where $u$ is the unit of the adjoint pair $(F,H)$
and $a$ is the action of $\Gamma$ on $H$.
For any $M\in\mod_{UB}(\ispk)$,
this map induces a natural map
\[ U{\rm Map}_c(G,M)\longleftarrow
   U{\rm Map}_c(G,UB\wedge_AM)\longleftarrow
   UB\wedge_A M\]
in $\mod_{UB}(\ispk)$.

\begin{lemma}
\label{lemma:endomorphism-monad}
If $\Psi_{UB}$ is a $k$-local equivalence,
then the composite map
$\Gamma\to  \Gamma HF\to HF$ is an equivalence
of functors.
\end{lemma}

\proof
It suffices to show that 
the induced map
$UB\wedge_AM\longrightarrow U{\rm Map}_c(G,M)$
is an equivalence 
for any $M$ in $\mod_{UB}(\ispk)$.
This follows from
Lemmas~\ref{lemma:bar-construction-representative-UB-A-M}
and \ref{lemma:psi-ub-imply-psi-m}.
\qqq

If $\Psi_{UB}$ is $k$-local equivalence,
by Lemma~\ref{lemma:endomorphism-monad}
and \cite[4.7.4]{Lurie2},
we see that 
the monad $\Gamma\in {\rm Alg}({\rm End}(\mathcal{C}))$ 
together with the object $\overline{U}H'\in 
{\rm Mod}_{\Gamma}({\rm Fun}(\mathcal{D},\mathcal{C}))$
is an endomorphism monad of $H$.
Hence we obtain an equivalence of quasi-categories
\[ {\rm Mod}_{\Theta}({\rm Fun}(\mathcal{D},\mathcal{C}))
   \stackrel{\simeq}{\longrightarrow} 
   {\rm Mod}_{\Gamma}({\rm Fun}(\mathcal{D},\mathcal{C}))  \]
compatible with the forgetful functors to 
${\rm Fun}(\mathcal{D},\mathcal{C})$,
and the object $\overline{H}\in 
{\rm Mod}_{\Theta}({\rm Fun}(\mathcal{D},\mathcal{C}))$ 
corresponds to
the object $\overline{U}H'\in 
{\rm Mod}_{\Gamma}({\rm Fun}(\mathcal{D},\mathcal{C}))$ 
under this equivalence.
Since the pair $(\Gamma,\overline{U}H')$ is an endomorphism 
monad of $H$,
in particular,
there is an equivalence
$\Gamma\stackrel{\simeq}{\to}\Theta$ in 
${\rm Alg}({\rm End}(\mathcal{C}))$.
This equivalence of algebra objects induces
an equivalence
\[ {\rm Mod}_{\Theta}(\mathcal{E})\stackrel{\simeq}{\longrightarrow}
   {\rm Mod}_{\Gamma}(\mathcal{E})\]
for any quasi-category left-tensored over ${\rm End}(\mathcal{C})$.
Taking $\mathcal{C}$ as $\mathcal{E}$,
we obtain the following theorem.
 
\begin{theorem}
\label{theorem:equivalence-Gamma-Theta}
If $\Psi_{UB}$ is $k$-local equivalence,
then there is an equivalence of quasi-categories
\[ \comod_{(UB,\Gamma)}(\ispk)\simeq 
   \comod_{(UB,\Theta)}(\ispk).\]
\end{theorem}

By Theorems~\ref{theorem:equivalence-BG-Gamma} and 
\ref{theorem:equivalence-Gamma-Theta},
we obtain the following corollaries.

\begin{corollary}
\label{corollary:mod-G-equivalent-comodules}
Let $G$ be a profinite group that has 
finite virtual cohomological dimension.
We assume that the localization functor
$L_k$ satisfies 
Assumption~\ref{condition:Behrens-Davis}.
If $\Psi_{UB}$ is a $k$-local equivalence,
then there is an equivalence of quasi-categories
\[ \mod_B(\isp(G)_k)\simeq
   \comod_{(UB,\Theta)}(\ispk).\]
\end{corollary}

Now we would like to show that 
the two formulations
of embeddings of module objects
in terms of model categories and
in terms of quasi-categories are equivalent
under some conditions.

First, we shall compare
the map $\mathcal{D}\to {\rm Mod}_{\Theta}(\mathcal{C})$
given by $\overline{H}\in 
{\rm Mod}_{\Theta}({\rm Fun}(\mathcal{D},\mathcal{C}))$
with the map
$\mathcal{C}(G)\to {\rm Mod}_{\Gamma}(\mathcal{C})$
given by $\overline{U}\in
{\rm Mod}_{\Gamma}({\rm Fun}(\mathcal{C}(G),\mathcal{C}))$.
The evaluation functor
$\mathcal{E}\times {\rm Fun}(\mathcal{E},\mathcal{C})\to
\mathcal{C}$ is a map of 
quasi-categories left-tensored over ${\rm End}(\mathcal{C})$
for any quasi-category $\mathcal{E}$.
This induces a map
\[ \mathcal{E}\times {\rm Mod}_T({\rm Fun}(\mathcal{E},\mathcal{C}))
\simeq
{\rm Mod}_T(\mathcal{E}\times {\rm Fun}(\mathcal{E},\mathcal{C}))
\longrightarrow
{\rm Mod}_T(\mathcal{C}) \]
for any monad $T\in {\rm Alg}({\rm End}(\mathcal{C}))$.
By adjunction, we obtain a map
\[ d(\mathcal{E},T):{\rm Mod}_T({\rm Fun}(\mathcal{E},\mathcal{C}))
   \longrightarrow
   {\rm Fun}(\mathcal{E},{\rm Mod}_T(\mathcal{C})),\]
which is an equivalence for any quasi-category $\mathcal{E}$.

We assume that $\Psi_{UB}$ is a $k$-local equivalence.
In particular, we have an equivalence 
$\Gamma\stackrel{\simeq}{\to}\Theta$
in ${\rm Alg}({\rm End}(\mathcal{C}))$.
By the naturality of the construction,
we obtain a commutative diagram
\begin{align}\label{diagram:LMod-Fun-commutation}
   \begin{array}{ccc}
    {\rm Mod}_{\Theta}({\rm Fun}(\mathcal{D},\mathcal{C}))&
    \stackrel{\mbox{$\scriptstyle d(\mathcal{D},\Theta)$}}
    {\hbox to 15mm{\rightarrowfill}}&
    {\rm Fun}(\mathcal{D},{\rm Mod}_{\Theta}(\mathcal{C}))\\
    \bigg\downarrow &  & \bigg\downarrow \\[3mm]
    {\rm Mod}_{\Gamma}({\rm Fun}(\mathcal{D},\mathcal{C}))&
    \stackrel{\mbox{$\scriptstyle d(\mathcal{D},\Gamma)$}}
    {\hbox to 15mm{\rightarrowfill}}&
    {\rm Fun}(\mathcal{D},{\rm Mod}_{\Gamma}(\mathcal{C})),\\
   \end{array}
\end{align} 
where all the arrows are equivalences.
We denote by $\widetilde{H}\in 
{\rm Fun}(\mathcal{D},{\rm Mod}_{\Theta}(\mathcal{C}))$
the image of $\overline{H}\in
{\rm Mod}_{\Theta}({\rm Fun}(\mathcal{D},\mathcal{C}))$
under the map
$d(\mathcal{D},\Theta)$
and 
by $\widetilde{U}\in{\rm Fun}(\mathcal{C}(G),
{\rm Mod}_{\Gamma}(\mathcal{C}))$
the image of $\overline{U}\in {\rm Mod}_{\Gamma}(
{\rm Fun}(\mathcal{C}(G),\mathcal{C}))$
under the map
$d(\mathcal{C}(G),\Gamma)$.
Since $\overline{H}$
corresponds to $\overline{U}H'$
under the equivalence
${\rm Mod}_{\Theta}({\rm Fun}(\mathcal{D},\mathcal{C}))
\stackrel{\simeq}{\to}
{\rm Mod}_{\Gamma}({\rm Fun}(\mathcal{D},\mathcal{C}))$,
we see that
$\widetilde{H}$
corresponds to $\widetilde{U}H'$
under the equivalence
${\rm Fun}(\mathcal{D},{\rm Mod}_{\Theta}(\mathcal{C}))
\stackrel{\simeq}{\to}
{\rm Fun}(\mathcal{D},{\rm Mod}_{\Gamma}(\mathcal{C}))$
using the commutative diagram~(\ref{diagram:LMod-Fun-commutation}).

The functor
${\rm Ex}=B\wedge_A(-): \mod_A(\spk)\to \mod_B(\sp(G)_k)$ 
of $\sp$-model categories
induces a functor 
\[ {\mathcal Ex}: \mod_A(\ispk)\to\mod_B(\isp(G)_k)\]
of quasi-categories.
We can identify the 
functor ${\mathcal Ex}^{\rm op}: \mod_A(\ispk)^{\rm op}\to
\mod_B(\isp(G)_k)^{\rm op}$
induced on the opposite quasi-categories
with $H':\mathcal{D}\to\mathcal{C}(G)$.
Recall that 
we have the functor
\[ {\rm Coex}: \mod_A(\ispk)\to \comod_{(UB,\Theta)}(\ispk) \]
of quasi-categories.
We can identify the functor
${\rm Coex}^{\rm op}:\mod_A(\ispk)^{\rm op}\to
\comod_{(UB,\Theta)}(\ispk)^{\rm op}$ 
induced on the opposite quasi-categories
with $\widetilde{H}:
\mathcal{D}\to {\rm Mod}_{\Theta}(\mathcal{C})$.
Furthermore, we recall that we have the map
$\mod_B(\isp(G)_k)\to 
\comod_{(UB,\Gamma)}(\ispk)$,
which is an equivalence under the assumptions
of Theorem~\ref{theorem:equivalence-BG-Gamma}.
We can identify the opposite of this map with 
$\widetilde{U}:\mathcal{C}(G)\to {\rm Mod}_{\Gamma}(\mathcal{C})$.
Since $\widetilde{H}$ corresponds to
$\widetilde{U}H'$ under the equivalence
${\rm Fun}(\mathcal{D},{\rm Mod}_{\Theta}(\mathcal{C}))
\stackrel{\simeq}{\to}
{\rm Fun}(\mathcal{D},{\rm Mod}_{\Gamma}(\mathcal{C}))$,
we obtain the following corollary.

\begin{corollary}
\label{corollary:equivalence-formulations}
Let $G$ be a profinite group that has 
finite virtual cohomological dimension.
We assume that the localization functor
$L_k$ satisfies 
Assumption~\ref{condition:Behrens-Davis}.
If $\Psi_{UB}$ is a $k$-local equivalence,
then there is an equivalence of functors
\[ {\mathcal Ex}\simeq {\rm Coex}\]
under the equivalence 
$\mod_B(\isp(G)_k)\simeq 
   \comod_{(UB,\Theta)}(\ispk)$.
\end{corollary}

This corollary shows that
the formulation of embeddings 
of module categories
in terms of model categories 
in Proposition~\ref{prop:fundamental-proposition-module-embedding} 
and that in terms of quasi-categories  
in Proposition~\ref{prop:quasi-category-module-embedding}
are equivalent.

\section{Embeddings over profinite Galois extensions}
\label{section:Galois-descent}

We assume that
$G$ is a profinite group which has finite virtual cohomological 
dimension.
Furthermore,
we assume that the localization functor
$L_k$ satisfies 
Assumption~\ref{condition:Behrens-Davis}.
In this section we show that
the two formulations of embeddings
of module categories are equivalent
if $\varphi: A\to B$ is a 
$k$-local $G$-Galois extension.

First, we show that a $k$-local $G$-Galois extension
gives an embedding of module categories.
Let $\varphi: A\to B$ be a $k$-local $G$-Galois extension.
We have a symmetric monoidal $\sp$-Quillen adjunction
\[ {\rm Ex}:\mod_A(\spk)\rightleftarrows\mod_B(\sp(G)_k): {\rm Re}\] 
by Lemma~\ref{lemma:symmetric-monoidal-module-adjunction-lemma},
which induces an adjunction of symmetric monoidal 
${\rm Ho}(\sp)$-algebras
\[ \mathbb{L}{\rm Ex}:
   {\rm Ho}(\mod_A(\spk))\rightleftarrows
   {\rm Ho}(\mod_B(\sp(G)_k)): \mathbb{R}{\rm Re}.\]
Let $\mathbf{T}$ be the full subcategory 
of $\Ho(\mod_A(\spk)$ consisting of $X$
such that the unit map 
$X\to \mathbb{R}{\rm Re}\mathbb{L}{\rm Ex}(X)$
is an isomorphism
\[ \mathbf{T}=\{X\in\Ho(\mod_A(\spk))|\ 
   X\stackrel{\cong}{\longrightarrow}
   \mathbb{R}{\rm Re}\mathbb{L}{\rm Ex}(X)\}.\]

\begin{proposition}
\label{prop:consistent-k-local-G-Galois-embedding}
If $\varphi: A\to B$ is a 
$k$-local $G$-Galois extension,
then the restriction of the functor
\[ \mathbb{L}{\rm Ex}:
   {\rm Ho}(\mod_A(\sp_k))\longrightarrow
   {\rm Ho}(\mod_B(\sp(G)_k)) \] 
to the full subcategory $\mathbf{T}$
is fully faithful as an $\Ho(\sp)$-enriched functor.
Furthermore,
if $\varphi$ is a consistent $k$-local $G$-Galois extension,
then the full subcategory $\mathbf{T}$
contains all dualizable objects.
\end{proposition}

\proof
The first part follows from
Propositions~\ref{prop:fundamental-proposition-module-embedding}.
If $B$ is a consistent $k$-local $G$-Galois extension
of $A$,
then $A\to B^{hG}$ is an equivalence
by \cite[Prop.~6.1.7(3) and 6.3.1]{Behrens-Davis}.
Hence the second part follows from
Proposition~\ref{prop:T-homotopy-fixedpoint-condition}.
\qqq

In the following of this section
we shall show that 
the underlying quasi-category of the model category
$\mod_B(\sp(G)_k)$
is equivalent to $\comod_{(UB,\Theta)}(\ispk)$.
Now we recall the construction of the map 
\[ \Psi_M: B(UB,A,M)\to U{\rm Map}_c(G,M)\]
for $M\in \mod_{UB}(\sp)$.
The map $\Psi_M$ is obtained 
by applying $U$ to the map
$B(B,A,M)\to {\rm Map}_c(G,M)$
in $\sp(G)$
that is adjoint to 
the map $UB(B,A,M)\cong B(UB,A,M)\to M$
induced by the action of $UB$ on $M$.

\begin{lemma}
\label{lemma:G-Galois-PhiM-eq}
If $\varphi: A\to B$ is a 
$k$-local $G$-Galois extension,
then the map 
$\Psi_M: B(UB,A,M)\to U{\rm Map}_c(G,M)$
is a $k$-local equivalence
for any cofibrant and fibrant object $M$
in $\mod_{UB}(\spk)$.
\end{lemma}

\proof
By the definition of $k$-local $G$-Galois extensions
\cite[Def.~6.2.1]{Behrens-Davis},
we have a
fundamental neighborhood system $\{U_{\alpha}\}$
of the identity element of $G$ consisting
of open normal subgroups and a directed system of 
finite $k$-local $G_{\alpha}$-Galois extensions  
$B_{\alpha}$ of $A$,
where $G_{\alpha}=G/U_{\alpha}$.
By the definition of finite Galois extensions
\cite[Def.~1.0.1]{Behrens-Davis},
we have a $k$-local equivalence 
\[ B_{\alpha}\wedge_A B_{\alpha}
   \stackrel{\simeq_k}{\longrightarrow} 
   {\rm Map}(G_{\alpha},B_{\alpha}).\]
Furthermore,
we have an isomorphism $B(B_{\alpha},A,B_{\alpha})\cong
B(B_{\alpha},A,A)\wedge_A B_{\alpha}$
and an equivalence 
$B(B_{\alpha},A,A)\stackrel{\simeq}{\to} B_{\alpha}$.
Since $A\to B_{\alpha}$ is a cofibration in 
the category of commutative symmetric ring spectra,
we obtain an equivalence
$B(B_{\alpha},A,B_{\alpha})\stackrel{\simeq}{\to}
 B_{\alpha}\wedge_A B_{\alpha}$
by \cite[Prop.~15.12]{MMSS}.
Hence we obtain a $k$-local equivalence
\[ B(B_{\alpha},A,B_{\alpha})
   \stackrel{\simeq_k}{\longrightarrow}
   {\rm Map}(G_{\alpha},B_{\alpha}).\]  

Let $r_{\alpha}: Q_{\alpha}M\to M$ be a 
cofibrant replacement in $\mod_{B_{\alpha}}(\spk)$
such that $r_{\alpha}$ is a trivial fibration.
We obtain a $k$-local equivalence
\[ B(B_{\alpha},A,Q_{\alpha}M)
   \stackrel{\simeq_k}{\longrightarrow}
   {\rm Map}(G_{\alpha},Q_{\alpha}M)\]
as in Lemma~\ref{lemma:psi-ub-imply-psi-m}.
Since $A$ and $B_{\alpha}$ are cofibrant
commutative symmetric ring spectra,
we see that $r_{\alpha}$
induces a $k$-local equivalence between
$B(B_{\alpha},A,Q_{\alpha}M)$
and $B(B_{\alpha},A,M)$
by using \cite[Prop.~15.12]{MMSS}.
Since $r_{\alpha}$ is a trivial fibration,
${\rm Map}(G_{\alpha},Q_{\alpha}M)\to
{\rm Map}(G_{\alpha},M)$
is also a trivial fibration.
Hence 
\[ B(B_{\alpha},A,M)\longrightarrow
   {\rm Map}(G_{\alpha},M) \]
is a $k$-local equivalence.
Since $\Psi_M$
is the colimit of the above maps
over the directed system,
the lemma follows from 
Proposition~\ref{prop:filteredcolimit-preserves-equivalences}.
\qqq

\begin{theorem}
\label{thm:Galis-mod-G-equivalent-comodules}
If $\varphi:A\to B$ is a $k$-local $G$-Galois extension,
then there is an equivalence of quasi-categories
\[ \mod_B(\isp(G)_k)\simeq
   \comod_{(UB,\Theta)}(\ispk).\]
Under this equivalence,
there is an equivalence of functors
\[ {\mathcal Ex}\simeq {\rm Coex}.\]
\end{theorem}

\proof
By Theorem~\ref{theorem:equivalence-BG-Gamma},
we have an equivalence
between $\mod_B(\isp(G)_k)$ and 
$\comod_{(UB,\Gamma)}(\ispk)$.
We can show that 
$\comod_{(UB,\Gamma)}(\ispk)$
is equivalent to $\comod_{(UB,\Theta)}(\ispk)$
as in Theorem~\ref{theorem:equivalence-Gamma-Theta}
by using Lemma~\ref{lemma:G-Galois-PhiM-eq}.
This completes the proof.
\qqq

Theorem~\ref{thm:Galis-mod-G-equivalent-comodules} 
shows that 
the two formulations of embeddings
of module categories are equivalent
if $\varphi: A\to B$ is a $k$-local $G$-Galois extension.




\end{document}